\documentclass[12pt,a4paper]{article}
\usepackage{amsmath, amscd, amssymb, amsthm, amstext, mathrsfs, xcolor}
\input{xypic}

\usepackage[final]{showkeys}

\theoremstyle{plain}
   \newtheorem{theorem}{Theorem}[section]
   \newtheorem{proposition}[theorem]{Proposition}
   \newtheorem{lemma}[theorem]{Lemma}
   
   \theoremstyle{definition}
   
   \newtheorem{definition}[theorem]{Definition}
   \newtheorem{example}[theorem]{Example}
   \theoremstyle{remark}
   
   \newtheorem{remark}[theorem]{Remark}
   \newtheoremstyle{citing}{3pt}{3pt}{\itshape}{\parindent}
{\bfseries}{.}{ }{\thmnote{#3}}
   \theoremstyle{citing} 
   \newtheorem*{varthm}{} 

\newcommand{\MM}{{M}}
\newcommand{\GG}{{\mathbf G}}
\newcommand{\KK}{{\mathbf K}}

\newcommand{\HH}{{\mathbf H}}

\newcommand{\EE}{{\mathbf E}}

\newcommand{\ZZ}{{\mathbb Z}}
\newcommand{\RR}{{\mathbb R}}
\newcommand{\CC}{{\mathbb C}}
\newcommand{\HHQ}{{\mathbb H}}

\newcommand{\Th}{{\mathop{Th}}}
\newcommand{\Bott}{{\mathop{BS}}}
\newcommand{\RP}{{\mathbb{R} \mathrm{P}}}
\newcommand{\CP}{{\mathbb{C} \mathrm{P}}}
\newcommand{\HP}{{\mathbb{H} \mathrm{P}}}
\newcommand{\CaP}{{\mathbb{O} \mathrm{P}^2}}

\newcommand{\geo}{{\gamma}}
\newcommand{\kgcycle}[2]{{\Gamma ({#2},{#1})}}
\newcommand{\kcycle}[1]{{\Gamma ({#1})}}

\newcommand{\norbd}{{\xi}}
\newcommand{\odb}[1]{{D^o ({#1})}}
\newcommand{\zun}{{\mathbf D}}

\newcommand{\Geo}{{\mathop{Geo}}}
\date{August 1, 2018}

\begin{document} 

\title{The suspended free loop space\\ of a symmetric space}
\author{Marcel B\" okstedt \& Iver Ottosen}
\maketitle

\begin{abstract}
Let $\MM$ be one of the projective spaces $\CP^n$, $\HP^n$ 
or the Cayley projective plane $\CaP$, and let
$\Lambda \MM$ denote the free loop space on $\MM$. Using 
Morse theory methods, we prove that the suspension spectrum 
of ${\Lambda \MM}_+$ is homotopy equivalent to the suspension 
spectrum of ${\MM}_+$ wedge a family of Thom spaces of explicit 
vector bundles over the tangent sphere bundle of $\MM$. 
\end{abstract}

\section{Introduction}
It has been known at least since Bott and Samelson \cite{BS} 
that it is possible to study the homotopy theory of the free 
loop space of a symmetric space using tools from differential
geometry and Lie group theory. One particularly
nice continuation of these ideas is the work by Ziller
\cite{Ziller}. He computes the integral homology of $\Lambda \MM$ 
for rank one symmetric spaces, and gives results and techniques
that apply to all symmetric spaces.  

In \cite{BO.SS} we approached the cohomology of a free
loop space $H^*(\Lambda \MM ;\ZZ /2)$ from a totally different 
angle. We used methods of cosimplicial spaces to set up
a spectral sequence which for simply connected
$\MM$ converges to this cohomology.
In order to compute the $E_2$-page of the spectral
sequence we needed a version of non-commutative 
homology, related to Andr\'e-Quillen homology of a ring. 
All of this is pure homotopy theory. There is nothing
in the methods that needs or uses that $\MM$ is a smooth 
manifold. The simplest non trivial case of the spectral 
sequence is when $M$ is a space such that 
$H^*(\MM ;\ZZ /2)$ is a truncated polynomial algebra. 
In \cite{BO} we used this approach to compute
$H^*(\Lambda \MM ;\ZZ /2)$ as a module over the
Steenrod algebra for such spaces. 

These calculations led to a strange observation.
In each case we considered, the cohomology splits as a sum
of finite dimensional modules. Moreover, we
recognized the pieces as cohomology of spaces
closely related to Thom spaces of iterated 
Whitney sums of the tangent bundle on $\MM$. 
We conjectured that the observed splitting was 
induced by an actual splitting of the suspension 
spectrum of $\Lambda \MM$ for these spaces.

The class of simply connected spaces such that
$H^*(\MM ;\ZZ/2)$ is a truncated polynomial algebra
is quite limited. The complex and quaternionic projective
spaces, together with the Cayley projective plane are not 
the only examples, but clearly stick out as the most 
interesting ones. And these are exactly the compact, 
simply connected, globally symmetric spaces of rank one.

The main results of this paper is that we can 
prove splittings of suspension spectra of
$\Lambda \MM_+$ for these rank one symmetric spaces. 
These theorems are given in section \ref{sec:Rank1}. 
For example:

\begin{varthm}
[Theorem \ref{th.complex.splitting}]
Let $p: S(\tau )\to \CP^n$ be the unit sphere bundle of the
tangent bundle $\tau$ on $\CP^n$. 
Let $\xi_m$ be the vector bundle 
$(p^*(\tau ))^{\oplus (m-1)} \oplus \epsilon$
over $S(\tau )$, where $\epsilon$ is a one dimensional trivial bundle.
Then, there is a homotopy equivalence of spaces
\[
\Sigma (\Lambda \CP^n)_+ \simeq 
\Sigma \CP^n_+ \vee \bigvee_{m=1}^\infty \Sigma \Th (\xi_m). 
\]
\end{varthm}

The proof is based on Ziller's methods. 
He does not quite formulate it in this way, but
he essentially proves that for a general globally
symmetric space, there is a splitting of 
the homology $H_*(\Lambda \MM ;\ZZ /2)$.
The main part of this paper consists of 
enhancing this argument to the point where we can prove
our conjecture on splittings of the suspension spectra.
In the tradition of Bott and Samelson, Ziller
used a mixture of differential geometry and
Morse theory. The extra ingredient we have added 
to this brew is a pinch of homotopy theory. 

The result we obtain is clearly stronger than 
Ziller's splitting of homology groups. For instance,
it follows that there are also splittings
for other homology theories. It also follows
that the splitting is compatible with the
action of cohomology operations. 

Previously, various splitting results have been
proved, for instance for the free loop space
of a suspended space. It is not clear to us that
this splitting is related to our splitting. 
A splitting for real projective space, which definitely is related
to ours, is shown in \cite{BCS}. 
   
Our proof is not completely formal. At one point,
we need a construction, which can only be
performed if a certain obstruction vanishes.
This obstruction lies in the cokernel of 
the map of representation rings induced
by a certain inclusion of Lie groups.
By calculation, this obstruction vanishes
in the cases we consider.

In an attempt to extend the splitting result to other symmetric spaces, we have later examined the 
infinite family of complex Grassmann manifolds. Here the methods of this paper only provide a stable splitting 
in certain special cases. It is hence unclear whether a stable splitting of $\Lambda \MM$ exists in general 
for compact symmetric spaces, but can not be obtained from the techniques presented here, or whether one only 
has a general homology splitting of $\Lambda \MM$ for some homology theories including singular homology.

In Section \ref{sec:Technical} we collect some background
material on equivariant differential topology and symmetric
spaces.

In section \ref{sec:Morse} we recall how Morse theory is applied
to free loop spaces. There is little original here.

In section \ref{sec:Bott} we reformulate and extend the 
method Bott and Samelson introduced for 
using the action of the isometry group to
construct interesting subspaces of the
free loop space. These spaces can be
identified as homogeneous spaces.
Much of this must have been known to Ziller.

Section \ref{sec:fixpt} contains the main new toy. We give 
methods to equivariantly split off the 
top cell of the homogeneous spaces
constructed in section \ref{sec:Bott}. This splitting is
often only possible in the stable situation, that
is after passing to the suspension spectrum. In
all cases, we can only construct the splittings
if certain obstructions vanish.
Then we show how we can use such splittings 
to obtain stable splittings of free loop spaces.

In section \ref{sec:Rank1}, we apply the theory of the
previous sections to the special cases
of projective spaces and to the Cayley projective plane.
We show that for these spaces the conditions
of section \ref{sec:fixpt} are satisfied, and we obtain splittings.
We also identify the summands of the splittings.
They are essentially Thom spaces of bundles 
over the total space of the unit sphere tangent
bundle over $\MM$.
 
Finally, in section \ref{sec:Thom} we check that the splitting we obtain here
agree with splitting conjectured in \cite{BO}.  

It is a pleasure to thank J. Tornehave for many discussions,
which have helped the paper substantially.

\section{Technical recollections}
\label{sec:Technical}

In this section we collect some background material on equivariant differential geometry and
symmetric spaces.

\subsection{Equivariant differential geometry}
We list some basic facts about actions of  
compact Lie groups on smooth manifolds. The tangent bundle of
a smooth manifold $\MM$ is denoted $\tau (\MM )$.

\begin{theorem}
\label{thm:sm.quotient}
Let $\GG$ be a compact Lie group and $\MM$ a smooth manifold 
equipped with a smooth and free right $\GG$-action. Then the 
orbit space $\MM/\GG$ has a unique smooth manifold structure, 
such that the canonical projection $\MM \to \MM /\GG$ is a submersion.
\end{theorem}

\begin{proof}
By the quotient manifold theorem \cite{Lee} 7.10, page 153 it suffices
to show that the action is proper. This follows directly since $\GG$ is compact
\cite{Lee} Corollary 7.2, page 147.
\end{proof}

Note that the real dimension of the quotient manifold is given by
\[\dim (\MM / \GG) = \dim (\MM ) - \dim (\GG ) . \]

\begin{theorem}
\label{thm:tb.borel}
Let $\GG$ be a compact Lie group with a closed subgroup $\HH$.
Assume that $\MM$ is a smooth manifold equipped with a smooth left $\HH$-action. Then 
the canonical projection 
\[ \pi : \GG \times_\HH \MM \to \GG /\HH ; \quad \pi ([g,x])=g\HH \]
is a morphism of smooth left $G$-manifolds. Furthermore, there is an isomorphism 
of $\GG$-vector bundles
\[ \tau (\GG \times_\HH \MM ) \cong 
\pi^* (\tau (\GG /\HH ))\oplus (\GG \times_{\HH } \tau (M)). \]
\end{theorem}

\begin{proof}
We have a smooth and free right action $\GG \times \MM \times \HH \to \GG \times \MM$ 
defined by $((g,p),h)\mapsto (gh,h^{-1}p)$. So $\GG \times_{\HH} \MM$ is a 
smooth manifold and the canonical projection $q: \GG \times \MM \to \GG \times_{\HH} \MM$ 
is a submersion by Theorem \ref{thm:sm.quotient}.
Similarly (when $M$ is a point) we have a smooth and free right action 
$\GG \times \HH \to \GG$ defined by $(g,h)\mapsto gh$ such that $\GG / \HH$ 
is a smooth manifold and the canonical projection $q^\prime : \GG \to \GG / \HH$ 
is a submersion.

The projection on the first factor $pr_1 : \GG \times \MM \to \GG$ is smooth, 
$\GG$-equivariant and there is a commutative diagram
\[
\xymatrix@C=1.0 cm{
\GG \times \MM  \ar[r]^-{pr_1} \ar[d]^-{q} & \GG \ar[d]^{q^\prime} \\
\GG \times_{\HH} \MM  \ar[r]^-{\pi} & \GG / \HH
}\] 
In particular, $q^\prime \circ pr_1$ is constant on the fibers of $q$ and since $q$ is a 
surjective submersion, \cite{Lee} Proposition 5.20, page 112 gives us that $\pi$ is smooth.

There are left $\GG$ actions 
$\GG \times \GG \times_\HH \MM \to \GG\times_\HH \MM$; 
$(g, [g^\prime , p]) \to [gg^\prime , p]$ and
$\GG \times \GG /\HH \to \GG / \HH$;
$(g, g^\prime \HH) \to gg^\prime \HH$. The projection $\pi$ is equivariant with 
respect to these. The actions are smooth as one sees by applying \cite{Lee} Propotion 5.20 
on two commutative squares having canonical projections vertically and actions by an 
arbitrary element of $\GG$ horizontally. 

By a standard result $q^\prime : \GG \to \GG / \HH$ is the projection map of a 
principal $\HH$-bundle. So
\[ \xymatrix@C=1.0 cm{
\MM \ar[r]^-{i_g} & \GG \times_{\HH } \MM \ar[r]^-{\pi} & \GG /\HH } \] 
where $g\in \GG$ and $i_g (p)=[g,p]$, is a fiber bundle. Since a fiber bundle is locally trivial, 
there exists an open neighborhood $U$ of $g\HH$ in $\GG /\HH$ and a fiber preserving 
diffeomorphism $\pi^{-1} (U) \cong U\times \MM$. So locally our fiber bundle
is isomorphic to a product bundle $\MM \to U\times \MM \to U$ and therefore we have an exact 
sequence for all $g\in \GG$ and $p \in \MM$ as follows:
\[ \xymatrix@C=1.0 cm{
0 \ar[r] & T_p (\MM ) \ar[r]^-{(i_g)_*} 
& T_{[g,p]}(\GG \times_\HH \MM) \ar[r]^-{\pi_*}
& T_{g\HH} (\GG /\HH ) \ar[r] & 0. } \]
In particular $\pi$ is a submersion. 

Recall that the tangent bundle $\tau (N)$ of a smooth $\GG$-manifold $N$ is a 
$\GG$-vector bundle, where the action of an element $g\in \GG$ on the total space 
is via its induced map $g_*$ of tangent spaces.

There is a commutative diagram of $\GG$-vector bundles
\[
\xymatrix@C=1.0 cm{
\tau (\GG \times_{\HH} \MM)  \ar[r]^-{\pi_*} \ar[d] & \tau( \GG /\HH ) \ar[d] \\
\GG \times_{\HH} \MM  \ar[r]^-{\pi} & \GG / \HH
} \] 
and hence a morphism into the pullback 
$\phi : \tau (\GG \times_\HH \MM ) \to {\pi }^* (\tau (\GG /\HH ))$.
Note that $\phi$ is surjective since $\pi$ is a submersion. 
In particular $\phi$ has constant rank and we get a short exact
sequence of $\GG$-vector bundles
\[ \xymatrix@C=1.0 cm{
 0 \ar[r] & \xi \ar[r] & \tau (\GG \times_\HH \MM) \ar[r]^-{\phi} 
& {\pi}^* (\tau (\GG /\HH )) \ar[r] & 0 } \]
where $\xi = \ker (\phi )$. Note that the dimension of $\xi$ equals 
the dimension of $\MM$ by this sequence. 

Since $\GG$ is a compact Lie group, \cite{Bredon} VI, Theorem 2.1, page 304, gives us that any 
smooth $\GG$-vector bundle has an invariant inner product.
Thus $\xi \subseteq \tau (\GG \times_{\HH} \MM )$ has a complementary $\GG$-vector bundle and
the short exact sequence splits. So we obtain an isomorphism of $\GG$-vector bundles
\[
\tau (\GG \times_\HH \MM ) \cong 
{\pi}^*(\tau (\GG /\HH ))\oplus \xi 
\]
and we only have to identify $\xi$.
 
There is a map $\psi: \GG \times \tau (\MM ) \to \tau (\GG \times_\HH \MM )$ given by 
$\psi (g,v_p) = (i_g)_*(v_p)$. By the identity
$i_{gh} \circ (h^{-1} \cdot) =i_g$ for $h \in \HH$ 
we have $(i_{gh})_* \circ (h^{-1} \cdot )_* = (i_g)_*$ so 
$\psi$ factors through $\GG \times_\HH \tau (\MM )$ and we obtain
a map over $\GG \times_\HH \MM$ as follows:
\[ \overline \psi : \GG \times_\HH \tau (\MM ) \to 
\tau (\GG \times_\HH \MM ). \]
From the first short exact sequence above we see that the image of 
$\overline \psi$ is contained in $\xi$. The dimension of its domain vector bundle is 
$\dim (\MM )$ which equals the dimension of $\xi$. By the first short exact sequence 
we also see that $\overline \psi$ is injective and the result follows.
\end{proof}

Recall that for a continuous left action of a topological group $\GG$ on a space $X$, the
isotropy group at $x\in X$ is defined as $\GG_x = \{ g \in \GG | gx=x \}$. 
It is a closed subgroup of $\GG$.

\begin{theorem}
\label{thm:orbit}
Let $\MM$ be a smooth manifold which is equipped with a smooth left action of a 
compact Lie group $\GG$. Then for any point $p\in \MM$ the orbit 
$\GG p \subseteq \MM$ is a smooth submanifold of $\MM$, and the map 
\[ \alpha_p : \GG /\GG_p \to \GG p; \quad g\GG_p \mapsto gp \] 
is a $G$-equivariant diffeomorphism.
\end{theorem}

\begin{proof}
This is Proposition I.5.4  of \cite{tD}, page 39. 
\end{proof}

\begin{theorem}
\label{thm:slices}
Let $\MM$ be a smooth manifold. Assume that $\HH$ is a compact Lie group which acts smoothly on  
$\MM$ with a fixed point $x\in \MM$. Then the following statements hold: 
\begin{enumerate}
\item The tangent space $T_x M$ is an $\HH$-representation. There exists an open neighborhood $U_x$ of 
$x$ in $\MM$, such that $hU_x\subseteq U_x$ for all $h\in \HH$, and an 
$\HH$-equivariant diffeomorphism $\phi :T_x\MM \to U_x$ with $\phi (0) = x$.
\item If $\HH$ is a closed subgroup of another compact Lie 
group $\GG$, then we can consider the $\GG$-equivariant section
\[ s:\GG /\HH \to \GG \times_\HH \MM ; \quad s(g\HH )= [g,x] \]
of the projection map $\pi$. The image $s(\GG /\HH)$ is a smooth submanifold 
of $\GG \times_\HH \MM$ whose normal bundle is $\GG$-equivariantly isomorphic to the vector bundle 
$$\GG \times_\HH T_x\MM \to \GG /\HH .$$
\end{enumerate}
\end{theorem}

\begin{proof}
For the first part, see \cite{tD} I.5.6, especially the reference to Bochner. 

For the second part, first note that because of Theorem \ref{thm:sm.quotient},
the three quotient spaces mentioned in the statement are smooth manifolds. 
Then consider the left action 
$$\GG \times (\GG \times_\HH \MM ) \to \GG \times_\HH \MM ,$$
and the point $[e,x] \in \GG \times_\HH \MM$. 
The isotropy group of this point is $\GG_{[e,x]} =\HH$ and its orbit is $\GG [e,x]=s(\GG /\HH )$.
By Theorem \ref{thm:orbit} we see that $s$ restricts to a $\GG$-equivariant diffeomorphism
$\GG / \HH \to s(\GG / \HH )$ and that $s(\GG /\HH )$ is a smooth submanifold of 
$\GG \times_{\HH} \MM$. We identify the normal bundle $\nu$ of this submanifold.

The map $\phi$, from the first statement, induces a $\GG$-equivariant diffeomorphism
\[
\GG \times_\HH T_x \MM \to \GG \times_\HH U 
\]
and $\GG \times_{\HH} U$ is an open neighborhood of $s(\GG/\HH)$. 
The normal bundle only depends on an open neighborhood of the submanifold so $\nu$ is
the normal bundle of $s(\GG/\HH)$ in $\GG \times_{\HH} U$. But via the equivariant 
diffeomorphism above this bundle is isomoprhic to the normal bundle of the zero section 
of the basse space in the total space of the $\GG$-vector bundle 
$\GG \times_\HH T_x\MM \to \GG /\HH$. This normal bundle is the $\GG$-vector bundle itself
\end{proof}

Finally, let us recall the Myers-Steenrod theorem \cite{Myers-Steenrod}:
\begin{theorem}
Let $\MM$ be a connected Riemann manifold. Then the isometry group $I(\MM )$ has the 
structure of a Lie group which acts smoothly on $\MM$. If $\MM$ is compact then $I(\MM )$
is also compact.
\end{theorem}

\subsection{Symmetric spaces}
\label{Section:Cayley}

There are several good expositions of the theory  
of symmetric spaces, for example in \cite{Helgason}, or in  \cite{Wolf}. 
Parts of the theory are also included in \cite{Milnor}.

A symmetric space is a connected Riemann manifold $\MM$ 
such that for each $p\in \MM$ there is an isometry $I_p :\MM \to \MM$ which leaves $p$ fixed and
reverses geodesics through $p$, i.e. if $\gamma$ is a geodesics and $\gamma (0) = p$ then 
$I_p(\gamma (t)) = \gamma (-t)$ for all $t$, \cite{Milnor} \S 20, page 109. 
Symmetric spaces are called globally symmetric spaces or Riemann globally symmetric spaces 
by some authors.

It follows from the definition that $\MM$ is complete \cite{Milnor}. The 
identity component $\GG = I_0 (\MM )$ of the isometry group $I(\MM )$ acts transitively on $\MM$
and for any point $p\in \MM$ the isotropy group $\GG_p$ is a compact subgroup. So there is a diffeomorphism $\GG /\GG_p \cong \MM$ such that $\MM$ is a homogeneous space, 
\cite{Helgason} IV, Theorem 3.3, page 208.

In section \ref{sec:Morse} we are going to study the 
Morse theory of closed curves on symmetric spaces.
This will lead us to questions on closed geodesic
curves and conjugated points along such curves. So as a preliminary and motivating
step, we will say a little about this subject.

First we present a result in a more general setting. Let $\MM$ be a connected and 
complete Riemannian manifold and fix a point $p\in \MM$. Let $\GG = I_0(\MM )$ 
and assume that the isotropy group $\KK = \GG_p$ is compact. 

The group $\KK$ acts on the unit sphere in the tangent space at $p$ via its differential 
\[ \KK \times S(T_p \MM ) \to S(T_p \MM ); \quad (g,v)\mapsto g_* (v). \]
On the other hand we have the set $\Geo_p (M)$ of unit speed geodesics $\gamma$ on 
$\MM$ with $\gamma (0) = p$. This set also has a $\KK$-action
\[ \KK \times \Geo_p (\MM ) \to \Geo_p (\MM ); \quad (g, \gamma )\mapsto g\gamma \]
and there is a $\KK$-equivariant bijection
\[ \Geo_p (\MM ) \xrightarrow{\cong} S(T_p \MM ) ; \quad \gamma \mapsto \gamma^\prime (0). \]
The inverse is given by the exponential map.

Fix a unit vector $v\in T_p M$ and let $\HH$ denote the isotropy group $\HH = \KK_v$.
Note that $\HH$ is the group that fixes the corresponding geodesics point wise. One has the 
following result regarding conjugated points:

\begin{lemma}
\label{lemma:group.conjugation}
Let $\gamma \in \Geo_p (\MM )$ with $\gamma^\prime (0) = v$.
Define the group $K_t$, for each $t>0$, by 
$\KK_t = \{ g \in \KK | g \gamma (t) = \gamma (t) \}$. Then the following holds:
\begin{enumerate}
\item If $\dim (\KK_t ) >\dim (\HH )$, then $\gamma (t)$ is conjugated to
$\gamma (0)$ along $\gamma$.
\item If $\dim (\KK_t )=\dim (\HH )$ and $\KK$ acts transitively on $S(T_p \MM)$, 
then $\gamma (t)$ is not conjugated to $\gamma (0)$ along $\gamma$.
\end{enumerate}
\end{lemma}

\begin{proof} 
Let $\exp_p : T_p \MM \to \MM$ denote the exponential map. We have that
$\gamma (t) = \exp_p (tv)$. By \cite{Milnor} Theorem 18.1, page 98, the point 
$\gamma (t)$ is conjugated to $\gamma (0)$ along $\gamma$ if and only if 
$tv$ is a critical point for the exponential map, i.e. if and only if the differential
\[ D_{tv}(\exp_p ) : T_{tv} (T_p\MM ) \to T_{\gamma (t)} \MM \] is not injective.
 
Also recall Gauss lemma: Let $S_t (T_p \MM) \subseteq T_p \MM$ denotes the set of 
tangent vectors at $p$ of length $t$. Then the image of the differential
\[
D_{tv}(\exp_p |) : T_{tv}(S_t (T_p \MM)) \to T_{\gamma (t)} \MM
\]
is orthogonal to the tangent vector $\gamma^\prime (t)$.
We explain how the lemma follows from these two results.

To prove the first statement we only use the first result.
The group $\KK$ acts on $S_t(T_p\MM )$ with isotropy group $\HH$ at the vector $tv$. 
The subgroup $\KK_t$ acts via the inclusion 
$\KK_t \hookrightarrow \KK$. By its definition, $\KK_t$ is a closed subgroup of $\KK$ so it is 
a compact Lie group. We have $\HH \subseteq \KK_t$ since elements in 
$\HH$ fixes $\gamma$ pointwise. Thus the isotropy group for the action 
of $\KK_t$ at $tv$ is again $\HH$. The composite
\[
\KK_t /\HH \hookrightarrow S_t(T_p\MM) \xrightarrow{\exp_p|} \MM; 
\quad k\HH \mapsto k_*(tv) \mapsto k\gamma (t)=\gamma (t)
\]
is a constant map. So it induces the trivial map on tangent spaces. 
By Theorem \ref{thm:orbit}, the differential of the
map to the left is injective. It follows that if $\dim (\KK_t ) >\dim (\HH )$, 
we obtain a non-trivial element in the kernel of $D_{tv}(\exp_p|)$ and hence also in 
the kernel of $D_{tv}(\exp_p)$.

We now consider the second statement of the theorem. 
Assume that $\gamma (r)$ is conjugated to $\gamma(0)$ along $\gamma$. 
Then $D_{rv}(exp_p)$ has a non trivial kernel. Decompose a nontrivial element in this kernel
as $w+s(rv)$, where $w\in T_{rv}(S_r(T_x\MM ))$ and $s\in \RR$. 
By Gauss lemma this orthogonal decomposition is preserved by the map $D_{rv}(\exp_p)$ so both 
components are mapped to zero. It follows that $s=0$ such that $w\neq 0$. 
We conclude that the kernel of $D_{rv}(\exp_p |)$ is non-trivial.

If $\KK$ acts transitively on the unit tangent sphere it also acts transitively on $S_t(T_p \MM )$. 
By Theorem \ref{thm:orbit} the action induces a diffeomorphism $\KK /\HH \cong S_t(T_p \MM)$. 
Furthermore, $\KK \to \KK / \HH$ is a submersion by Theorem \ref{thm:sm.quotient}. So $w$ lifts
to an element $\widetilde{w}$ in the Lie algebra $T_e \KK$. Let $\phi : \RR \to \KK$ be the corresponding 
1-parameter group. 

Note that $\phi$ does not fix $\gamma$ pointwise. If it did, $\phi$ would be 
a 1-parameter group of $\HH$, and $\RR \xrightarrow{\phi} \KK \xrightarrow{Q} \KK /\HH$ a constant map.
But $Q_* (\widetilde{w} ) \neq 0$.

However, $\phi (s) \gamma (0) = \gamma (0)$ for all $s$ by definition of $\KK$. 
So it suffices to show that there exists an $\epsilon > 0$ such that  
$\phi (s) \gamma (r) = \gamma (r)$ for all $|s|<\epsilon$, since then $\dim(\KK_r) > \dim(\HH)$.
In fact, it suffices to show that
\begin{equation}
\label{zeroderivative}
\frac d {ds} (\phi (s) q) |_{s=0} = 0
\end{equation}
by the following argument:
Put $q=\gamma (r)$ and let $X$ denote the vector field on $\MM$ generated by $\phi$ (see \cite{Milnor} page 10).
Then the curve $\sigma (s) = \phi (s) q$ satisfies the differential equation ${\sigma }^{\prime } (s) = X_{\sigma (s)}$
with initial condition $\sigma (0) = q$  (see \cite{Milnor} page 10 proof of lemma 2.4). Under the above assumption 
on the derivative at zero, the constant curve with value $q$ is also a solution of this initial value problem. So
$\sigma (s) = q$ in a neighbourhood of zero by uniqueness of solutions. 

Regarding (\ref{zeroderivative}), note that the exponential map is $\KK$-equivariant, such that we have a 
commutative diagram
\[
\xymatrix@C=1.0 cm{
\RR \ar[r]^-{\phi} \ar@{=}[d] & \KK \ar[r]^-{i} \ar@{=}[d] & \KK \times T_p \MM \ar[r] \ar[d]^-{id\times exp_p} 
& T_p \MM \ar[d]^-{exp_p} \\
\RR \ar[r]^-{\phi} & \KK \ar[r]^-{j} & \KK \times \MM \ar[r] & \MM 
} \] 
where $i(k)=(k, rv)$ and $j(k)=(k, q)$. The top composite is the map $s\mapsto \phi (s) rv$ 
and the bottom composite $s\mapsto \phi (s) q$. Consider the induced diagram on tangent spaces with
$T_0 \RR$ in the upper left corner. In that diagram the upper composite maps $1$ to $w$
and the lower composite maps $1$ to the desired derivative. Since $w$ is in the kernel of $D_{rv}(exp_p)$
the result follows.
\end{proof}

\begin{remark}
In a symmetric space, the relation between the isometry
group and conjugated points is even better, since the
isometry group is big. We will see a version of this improved
correspondence in section \ref{sec:Bott} and \ref{sec:deformation}, 
especially Proposition \ref{prop:Bott-Samelson} 
and Theorem \ref{dimension.negative.bundle}. 
\end{remark}

In section \ref{sec:Rank1} we are going to study 
symmetric spaces of rank 1. We continue this paragraph 
by discussing these, especially the projective Cayley 
plane which is the most exotic one. A part of this is written in 
expository style, but we give proofs or precise references for the results we 
actually need in section \ref{sec:Rank1}. There does not seem to be a 
reference to the projective Cayley plane that covers it all, 
so we have to collect results from several sources.

\begin{definition}
\label{def:isotropic}
A Riemannian manifold $\MM$ is said to be \emph{isotropic}, if the
isometry group of $\MM$ acts transitively on the total space of the 
unit sphere bundle of the tangent bundle.
\end{definition}

The isotropy condition is quite constraining.
It is equivalent to the condition that $\MM$ is an Euclidean space or
a symmetric space of rank 1, and it is also equivalent to the
condition that the group $\GG = I_0 (\MM )$ acts transitively on
pairs of points with a fixed distance.
These equivalences are proved in \cite{Wolf} Corollary 8.12.9.

According to \cite{Wolf} 8.12.2 this class includes
spheres, complex and quaternionic projective spaces and 
the Cayley projective plane, but no other
compact, simply connected spaces.

Let $\MM$ be a compact symmetric space of rank 1. Since it is compact, 
it has a finite diameter. By compactness, there exists two points on $\MM$
with distance equal to the diameter. By the Hopf-Rinow 
theorem, these two points are joined by a minimal geodesic 
of length equal to the diameter of $\MM$. Define the {\em antipodal set} of a point 
$x\in \MM$ by 
\[ A_x=\{ y\in \MM | d(x,y)=\text{diam} (\MM ) \} .\]
According to \cite{Helgason} VII, \S 10, page 328, the isotropy 
group $\KK= \GG_x$ acts transitively on $A_x$. So by Theorem \ref{thm:orbit}, 
$A_x$ is a compact submanifold of $\MM$, diffeomorphic to
$\KK/\KK_1$, where $\KK_1 \subseteq \KK$ is the isotropy group of some point $y\in A_x$.

\begin{lemma}
\label{lemma:geodesics.symmetric}
Let $\gamma$ be a closed geodesics on $\MM$ parametrized by arc length $t$.
The length of $\gamma$ equals $2mL$ for some positive integer $m$. 
The point $\gamma (L)$ has maximal distance from 
$\gamma (0)$, so $L = \text{diam} (\MM )$ and $\gamma (L) \in A_{\gamma (0)}$.
Let $\KK = \GG_{\gamma (0)}$, $\KK_1 = \KK_{\gamma (L)}$ and let $\HH$ be the
subgroup that fixes $\gamma$ point wise.
There are at most three conjugacy classes of isotropy groups
of arbitrary pairs of points in $\MM$, represented by $\HH \subseteq \KK_1\subseteq \KK$, and
the following statements hold:
\begin{enumerate}
\item The point $\gamma (t)$ is not conjugated to $\gamma (0)$ along $\gamma$ 
if $t$ is not an integer multiple of $L$.
\item The point $\gamma (2kL)$ is conjugated to $\gamma(0)$ along
$\gamma$ for any integer $k$ and $\gamma (2kL)=\gamma (0)$.
\item The point $\gamma ((2k+1)L)$ is conjugated to $\gamma (0)$ along $\gamma$ 
for every integer $k$ if $\HH$ has positive codimension in $\KK_1$.
\end{enumerate}
\end{lemma}

\begin{proof} 
There exist a simple closed geodesics, and the 
simply closed geodesics are all of the same length by 
\cite{Helgason} VII, Proposition 10.2. Let $L$ be the half of that length.
Because $\GG$ acts transitively on geodesics, all closed geodesics will
run through a simple closed geodesics an integral number of times.
Since any two points are connected by a minimal geodesic, this also
shows that the diameter of $\MM$ is at most $L$.

The exponential map is a diffeomorphism from the open ball
in the tangent space 
$B_L(x)=\{v \in T_x\MM | |v| < L \}$ onto the open set
$\MM \setminus A_x$ according to \cite{Helgason} VII 
Theorem 10.3. In particular, this shows that the diameter of $\MM$ is exactly $L$. 
 Since the exponential map is $\KK$-equivariant, it follows that
the isotropy subgroup of $\KK$ fixing a point in
in $B_L(x)\setminus \{ 0 \}$ is exactly $\HH$.

Since $\KK$ acts transitively on the unit tangent
sphere at $x\in \MM$, lemma \ref{lemma:group.conjugation}
tells us that a point $\gamma (t)$ cannot be conjugated
to $\gamma(0)$ unless either $t=2kL$, such that $\gamma (t) = \gamma (0)$, 
or $t=(2k+1)L$, such that $\gamma (t) \in A_{\gamma (0)}$,
for some integer $k$. 

The isotropy group of $\gamma (2kL) = \gamma (0)$ is $\KK$ and $\dim (\KK ) > \dim (\HH )$
since $\dim (\KK)-\dim (\HH ) = \dim (\KK / \HH) = \dim (T_x\MM )-1$. So lemma \ref{lemma:group.conjugation} tells us that these points are conjugated to $\gamma (0)$.

The isotropy group of $\gamma ((2k+1)L) = \gamma (L)$ is $\KK_1$.
So by lemma \ref{lemma:group.conjugation}, these points are conjugated to $\gamma (0)$ 
if and only if $\dim (\KK_1 ) > \dim (\HH )$.

Finally, we show the statement regarding conjugacy classes. For a general group action 
$G \times A \to A$ one has $g G_a g^{-1} = G_{ga}$ for $g\in G$ and $a\in A$. 

Let $x$ and $y$ be a pair of points in $\MM$. If $x=y$ then there is an element 
$g\in \GG$ with $gx=\gamma (0)$ since $\GG$ acts transitively on $\MM$. 
So $\GG_x \cap \GG_y=\GG_x$ is conjugate to $\KK$.  
If $d(x,y)=s$, where $0<s\leq L$, then there is an element $g\in \GG$ such that $gx=\gamma (0)$
and $gy=\gamma (s)$ by the transitivity result from \cite{Wolf} mentioned above. 
Thus $\GG_x \cap \GG_y$ is conjugate to $\GG_{\gamma (0)} \cap \GG_{\gamma (s)}$.
If $s<L$ this group equals $\HH$ and if $s=L$ it equals $\KK_1$.
\end{proof}

\begin{remark}
Actually, we will always have that $\dim (\KK_1 )>\dim (\HH )$.
\end{remark}

We now specialize this to the Cayley projective plane $\CaP$.
The exotic Lie group $F_4$ contains $Spin(9)$ as a subgroup. 
The Cayley projective plane is the homogeneous quotient $\CaP \cong F_4/Spin(9)$.
is an symmetric space of rank 1. This is shown in  
\cite{Wolf} Theorem 8.12.2. The isotropy group $Spin(9)$ 
acts on the tangent space $T_x\CaP$. This is a real 16 dimensional representation.

\begin{remark}
We won't need to determine which representation of $Spin(9)$ 
this is. But one can show that it is the spinor representation $R_9$
of $Spin(9)$ as claimed in \cite{Wolf}, proof of Theorem 8.12.2.
\end{remark}

What we will need are the following facts.

\begin{theorem}
\label{th:Cayley.geodesics}
A minimal closed geodesic $\gamma$ on $\CaP$ has exactly two points
$\gamma (L)$ and $\gamma (2L)$ conjugated to $\gamma (0)$.
Let $\KK_2$ be the isotropy group of $\gamma(0)$, and 
$\KK_1$ be the subgroup of $\KK_2$ which fixes the point
$\gamma (L)$.
There are group isomorphism $\phi_1 :\KK_1 \to Spin(8)$
and $\phi_2 : \KK_2 \to Spin(9)$ so that the composite
\[
Spin(8) \xrightarrow{\phi_1} \KK_1 \subset \KK_2
\xrightarrow{\phi_2^{-1}} Spin(9)
\]
is the standard inclusion.
\end{theorem}

\begin{proof}
The calculations of isotropy groups is given
in \cite{BP} Chapter 1, \S 5, (3). The model for $\CaP$ used
in this article is Freudenthals model, see also \cite{Freudenthal}.  
Points in $\CaP$ are identified with projections in a Jordan algebra
of $3\times 3$ matrices over the octonions. They identify this space 
with the unique symmetric space of type $F_4/Spin(9)$.
See \cite{BP} Chapter 3, \S 10, at the end of the paragraph.

They pick three points $E_1,E_2,E_3$ with the property that 
$E_1+E_2+E_3=1$. They compute that the group fixing
$E_1$ is $Spin(9)$. The group fixing $E_1$ and $E_2$ actually 
fixes all three of them since $E_1+E_2+E_3=1$ is invariant under
the action of the isotropy group $F_4$. In Chapter 1, \S 5, (3) they 
compute this group to be $Spin(8)$, and that the inclusion map of 
this in $Spin(9)$ is the standard inclusion. The fact that this is 
the standard inclusion is important to us. This was certainly
also known to Freudenthal. A more detailed argument
for that infinitesimally the inclusion of root systems is the correct 
one is given in \cite{Freudenthal} 4.12. But this information 
determines the inclusion.

It follows from lemma \ref{lemma:geodesics.symmetric}
that there are at most two points conjugated to $\gamma (0)$.
Recall that $\dim (Spin (n)) = n(n-1)/2$.
Since $S(T_x(\CaP ))\cong Spin(9)/\HH$ we have that 
$\dim (\CaP )-1 = \dim (Spin(9))-\dim (\HH )$ such that
$\dim (\HH ) = 21$. Thus, 
$\dim (\KK_1 )=\dim (Spin(8))= 28>\dim (\HH )$, 
and by lemma \ref{lemma:geodesics.symmetric} we conclude that
both $\gamma (L)$ and $\gamma (2L)$ are conjugated to $\gamma (0)$.
\end{proof}

\begin{remark}
We are actually not going to need the following results,
so we state them without proof. The group $\HH$ is isomorphic
to $Spin(7)$. The inclusion $\HH \subset \KK_1$ 
is not standard. It is the composite of the usual inclusion with the
``triality'' automorphism of $Spin(8)$;
\[
Spin(7) \subset Spin(8) \xrightarrow{\theta} Spin(8).
\]
The inclusion $\HH \subset \KK_2$ is the composite
\[
Spin(7) \subset Spin(8) \xrightarrow{\theta} Spin(8)\subset Spin(9).
\]
The space $\KK_1 /\HH$ is the unit sphere in the
standard representation of $Spin(7)$, and the space 
$\KK_2/\HH$ is the unit sphere in the tangent representation of 
$Spin(7)$ at $[e]\in F_4/Spin(9)$. As we mentioned above, the 
$Spin(9)$ representation $\mathfrak{f}_4 /\mathfrak{spin} (9)$ is
the spin representation.
\end{remark}

\section{Morse theory and geodesics on symmetric spaces}
\label{sec:Morse}

Morse theory for free loop spaces is a variation on Morse theory for based
loop spaces, developed by Bott, Samelson and later Klingenberg.
The specialization of this theory to the particularly
agreeable case of symmetric spaces is mostly due to 
Bott, Samelson and Ziller. 

\subsection{Morse theory for free loop spaces}

There is a version of Morse theory for the free loop space of a compact Riemann manifold.
We use Klingenberg's book \cite{Klingenberg} as the standard reference. The parts which we 
refer to are not controversial. 

Let $\MM$ be a compact and connected Riemann manifold and let $S^1$ denote the 
circle $\RR /\ZZ$. The free loop space $\Lambda \MM $ is the space of 
maps $\gamma : S^1 \to \MM$ of Sobolev class $H^1$ i.e. 
absolutely continuous maps $\gamma :S^1\to M$, such that 
$\int_0^1 |\gamma^\prime (t)|^2 dt <\infty$. 
One can show that $\Lambda \MM $ has the structure of a smooth 
Hilbert manifold. There are inclusions for the smooth, piecewise smooth and continuous versions 
of the free loop space 
\[ C^\infty (S^1, \MM ) \subseteq C^{\infty}_{pw} (S^1, \MM ) \subseteq \Lambda \MM \subseteq C^0(S^1, \MM ) \]
and these inclusions are homotopy equivalences. 
For detalis, see \cite{Klingenberg} \S 1.1 and \S 1.2. 

The energy functional on the free loop space is defined by
\[ E: \Lambda \MM \to \RR ; \quad 
E(\gamma ) = \frac 1 2 \int_0^1 |\gamma^\prime (t)|^2 dt . \]
It is a smooth map. The critical points for the energy functional are precisely the closed 
geodesic curves on $\MM$, \cite{Klingenberg} \S 1.3.

In the simplest version of Morse theory, critical points are isolated. However,
there is a continuous action of the circle group $S^1$ on $\Lambda \MM$ 
given by rotation of loops and the energy functional is invariant under this action.
Every element in the orbit $S^1\gamma$ of a periodic geodesics $\gamma$ is a periodic 
geodesics of the same energy level. So a non-constant periodic geodesics cannot be an 
isolated critical point (constant ones are clearly not isolated as they lie on $\MM\subseteq \Lambda M$). 
Furthermore, we have a smooth action of the group $\GG =I_0(\MM )$ on 
$\Lambda \MM$. Again the orbit $\GG \gamma$ consists of periodic geodesics 
of the same energy level. So $\GG \gamma$ form a critical submanifold of 
$\Lambda \MM$ and we will be considering cases where this submanifold is positive 
dimensional. 

So assume that critical points form critical smooth submanifolds of $\Lambda \MM$. 
Let $N$ be one of these critical submanifolds. Also assume that we are in the special case where
the adjoint $S^1\times N \to \MM$ of the inclusion map $i: N\hookrightarrow \Lambda \MM$ possess 
a subdivision $0=t_0<t_1< \dots <t_k=1$ such that 
$[t_{j-1},t_j]\times N \to \MM$ is smooth for $j=1, 2, \dots ,k$.
That is, each loop is assumed to be piecewise differentiable.
All curves in the same connected component of the critical manifold 
will have the same energy since the energy integral is continuous and locally constant on $N$. 

The inclusion of $N$ induces a tangent map 
$i_*: T_\gamma (N)\to T_\gamma (\Lambda \MM )$.
The tangent vector space $T_\gamma (\Lambda \MM )$
can be considered as the space of $H^1$-vector fields along $\gamma$. 
Integration along $\gamma$ induces an inner product on 
$T_\gamma (\Lambda \MM )$, eventually induced from the metric on $\MM$.

There is a 'normal bundle' $\norbd$ defined on
each component of $N$. The fiber of $\norbd$ at 
$\gamma \in N$ is the vector space of periodic 
vector fields along $\gamma$.
 
Assume that the critical manifolds satisfy 
the Bott-Morse non-degeneracy condition. This says that the
null space of the Hessian of the energy functional is
exactly the tangent directions of $N$, that is, the image of $i_*$.
The metric on  $T_\gamma (\Lambda \MM )$ induces a splitting 
$\norbd \cong \norbd^- \oplus \norbd^+$ of vector bundles over $N$.
The Hessian of the energy function is
positive definite on $\norbd^+$, and negative definite on
$\norbd^-$. The bundle $\norbd^+$ is infinite dimensional,
but $\norbd^-$ is a finite dimensional vector bundle.

Let $\Lambda^a \MM$ be the subspace of loops of energy
less than or equal to $a$. Assume that that $k$ is a critical value and that there 
are no critical values except $k$ in the interval 
$[k-\epsilon,k+\epsilon]$ where $\epsilon >0$. 
Also assume that all the critical points in
$E^{-1}([k-\epsilon ,k+\epsilon ])$ are situated on 
a connected, non-degenerate critical submanifold $N$. 
This manifold will then be isolated.
 
The main statement of Morse theory of the free loop space is
\cite{Klingenberg} 2.4.10:
\begin{theorem}
\label{th.morse.theory}
There is a homotopy equivalence
\[
\Lambda^{k+\epsilon} \MM \simeq 
\Lambda^{k-\epsilon} \MM \cup_\sigma D(\norbd^-)
\]
for some gluing map 
$\sigma : S(\norbd^-) \to \Lambda^{k-\epsilon } \MM$.
\end{theorem}

Of course, one can handle the case of 
several isolated critical manifolds in the same fashion.
Assume that $k$ is the only critical value of $E$ 
in the interval $[k-\epsilon ,k+\epsilon ]$. 
Suppose that the critical point with critical value
$k$ is a union of non degenerate critical
manifolds $N_\nu$, $\nu = 1, 2, \dots , m$ each with a negative bundle
$\norbd^-_\nu$. Then we have a homotopy equivalence
\[
\Lambda^{k+\epsilon} \MM \simeq 
\Lambda^{k-\epsilon} \MM \cup_{\sigma_1} D(\norbd^-_1)\cup_{\sigma_2} \dots 
\cup_{\sigma_m} D(\norbd^-_m).
\]

Also assume that we are in the special case where
the adjoint $S^1\times N \to \MM$ of 
the inclusion map $i: N\hookrightarrow \Lambda \MM$ possess 
a subdivision $0=t_0<t_1< \dots <t_k=1$ such that 
$[t_{j-1},t_j]\times N \to \MM$ is smooth for $j=1, 2, \dots ,k$.
That is,  each loop is assumed to be piecewise differentiable.

We will need more precise information about the 
gluing map. Let $\eta$ be a vector bundle over $N$. We let
$\odb \eta$ denote the open unit disk bundle and
$s: N \to \odb \eta$ the $0$-section.
Assume that $f: \odb \eta \to \Lambda \MM$ is a continuous map,
such that $f\circ s=i$ and such that the adjoint
$S^1\times \odb \eta \to M$ is smooth in the $\odb \eta$ direction
and piecewise smooth in the $S^1$ direction (with respect to
some subdivision of the subdivision as we used for the inclusion map $i$).

\begin{lemma}
\label{lemma:negative.bundle}
Assume that the dimension of $\eta$ equals the dimension
of the negative bundle $\norbd^-$. Assume also that 
the Hessian of $E\circ f$ is negative definite on 
each fiber of $\eta$. Then, there is a homotopy equivalence
$$\Th (\eta)\simeq \Lambda^{k+\epsilon} \MM /\Lambda^{k-\epsilon} \MM .$$
\end{lemma}

\begin{proof}
The Hessian is functorial in the following sense: If $h: M_1\to M_2$ and $g:M_2 \to \RR$ are smooth maps
and $p\in M_1$ maps to a critical point $h(p)$ for $g$, i.e. $D_{h(p)}(g)=0$, then $p$ is a critical 
point for $g\circ h$ (by the chain rule) and
\[ H_p(g\circ h) (v,w) = H_{h(p)}(g) (D_p(h)(v), D_p(h)(w)) . \]
This result is a consequence of the description of the Hessian given in \cite{Milnor} page 74. 

Thus for any $\gamma \in N$ we see that $s(\gamma )$ is a critical point for $E\circ f$ and
\begin{equation} \label{Hessian nat}
H_{s(\gamma)} (E\circ f)(v,w) = H_{i(\gamma )}(E)(D_{s(\gamma)}(f)(v), D_{s(\gamma )}(f) (w)) .
\end{equation}
Our assumptions imply that that $N$ is a non-degenerate critical submanifold for the function $E\circ  f$. 
That is the null space of $H_{s(\gamma )}(E\circ f)$ equals the image of $D_\gamma (s)$. 
The inclusion $\supseteq$ follows by (\ref{Hessian nat}) since $N$ is a non-degenerate critical submanifold for $E$.
For the other inclusion $\subseteq$ we use the splitting 
\[ T_{s(\gamma )} (\odb \eta ) \cong T_\gamma (N) \oplus \eta_\gamma \]
where $\eta_\gamma$ denotes the fiber of $\eta$ over $\gamma$, and the assumption that the Hessian
is negative definite on $\eta_\gamma$.

Now consider the restriction  
$f| : {\odb \eta}_\gamma \to \Lambda \MM$.
The differential of this restriction is an injective map
$f|_* : \eta_\gamma \to T_\gamma \Lambda \MM$ by (\ref{Hessian nat}) and the negative definite assumption. 
These maps combine to an inclusion of vector bundles $\eta \subset \norbd$
which again by (\ref{Hessian nat}) restricts to an inclusion $\eta \subset \xi^-$. 
Since we are assuming that $\eta$ and $\norbd^-$ have the same 
dimension, this composite is an isomorphism of vector bundles $\eta \cong \xi^-$.

By the proof of 2.4.11. in \cite{Klingenberg} one has a homotopy equivalence 
\[ \Th (\norbd^-) \simeq 
\Lambda^{k+\epsilon} \MM / \Lambda^{k-\epsilon} \MM \]
and by the isomorphism above we have $\Th (\eta ) \cong \Th (\norbd^- )$.
\end{proof}

\subsection{The free loop space of a symmetric space}

From now on, let $\MM$ be a compact symmetric space. 
Due to the existence of a large isometry group, Morse 
theory has some very special properties in this situation.
Ziller \cite{Ziller} computes the critical submanifolds and shows 
that they satisfy the Bott-Morse non-degeneracy condition. 

Let $\sigma :[0,a] \to \MM$ be a 
geodesic parametrized by arc length which is also
a closed curve $\sigma (0)=\sigma (a)$. It turns out that the tangent vectors at the endpoint agree,
$\sigma^\prime (0)=\sigma^\prime (a)$, so that 
$\gamma$ is actually a periodic geodesic, \cite{Ziller} page 5.  

The identity component of the isometry group $\GG = I_0(\MM )$ acts
smoothly on $\Lambda \MM$. Let $\gamma$ be a periodic geodesics on $\MM$ and
let $\HH \subseteq \GG$ be the subgroup that fixes $\gamma$ pointwise. 
Then the orbit $\GG \gamma$ is a submanifold of $\Lambda \MM$ which 
is diffeomorphic to $\GG /\HH$ by Theorem \ref{thm:orbit}.
 
Furthermore, by \cite{Ziller} Theorem 2 one has the following result:

\begin{theorem}
Let $\MM$ be a symmetric space and let $\gamma$ be a periodic geodesics on $\MM$. 
Then $\GG \gamma$ is a non-degenerate critical submanifold of $\Lambda \MM$.
\end{theorem}

Obviously every critical point of $E$ is contained in one of these 
critical submanifolds. 
We say that two closed geodesics $\gamma$ and $\gamma^\prime$ on $\MM$ 
are $\GG$-equivalent if $g\gamma = \gamma^\prime$ for some $g\in \GG$. This is an equivalence
relation where the equivalence classses are the orbits which we have identified as the critical 
submanifolds. In particular, the compact critical submanifolds are disjoint and hence isolated. 
By Theorem \ref{th.morse.theory} it follows that $\Lambda \MM$ has a filtration, with filtration
quotients the Thom spaces of the negative bundles $\Th (\norbd^- )$.
The filtration is indexed by the value of the energy functional 
on the components of the critical sets. If there are several 
components of the same energy level, we chose an order of those, 
and reindex our filtration. The filtration quotients are exactly $\Th (\norbd^- )$.

This means that we have the following:
\begin{theorem}
\label{th.filtration1}
Let $\MM$ be a compact symmetric space. There is a filtration 
$$\dots \subseteq F^{\nu^\prime} \subseteq F^{\nu } \subseteq \dots \subseteq \Lambda \MM $$
indexed by the $\GG$-equivalence classes of geodesic loops. These equivalence classes are 
ordered in such a way that the energy is weakly increasing. For two consecutive steps 
$F^{\nu^\prime } \subseteq F^\nu$
in the filtration, there is a homotopy equivalence
$F^\nu /F^{\nu^\prime } \simeq \Th (\norbd_\nu^- )$
where $\norbd_\nu^-$ is the negative bundle over the critical submanifold indexed by $\nu$.
\end{theorem}

\section{The Bott-Samelson map}
\label{sec:Bott}

This paragraph is a reformulation of \cite{Ziller}, \S 3.1.
We intend to analyze closer how $F^\nu$ is built out of
$F^\nu / F^{\nu^\prime }$ and $\Th (\norbd^-_\nu )$.
The tool is the construction of a subspace of broken closed 
geodesics inside $F^\nu$. The construction is due to
Bott and Samelson \cite{BS}. 

\subsection{$K$-cycles}
The space of broken geodesics is related to certain homogeneous spaces called  
Bott-Samelson $K$-cycles. They are defined as follows \cite{BS}, Chapter I, 4.4:

\begin{definition}
\label{definition:K-cycle}
Let $\GG$ be a compact Lie group with a closed subgroup $\HH$.
Let $\KK_\bullet =(\KK_1, \dots ,\KK_m)$ be an $m$-tuple of 
closed subgroups such that $\HH \subset \KK_i \subset \GG$ for 
$1\leq i\leq m$. Define a right action of 
$\HH^m = \HH \times \dots \times \HH$ on 
$\KK_1 \times \dots \times \KK_m$ by
$$(c_1, \dots ,c_m)*(a_1, \dots ,a_m) = 
(c_1a_1, a_1^{-1}c_2a_2, a_2^{-1}c_3a_3, \dots ,
a_{m-1}^{-1}c_ma_m).$$
The $K$-cycle is the associated orbit space
$$\EE (\KK_\bullet ;\HH ) = (\KK_1 \times \dots \times \KK_m )/\HH^m.$$
\end{definition}

\begin{lemma}
\label{lemma:K-cycle}
The following properties hold:
\begin{enumerate}
\item The space $\EE (\KK_\bullet ;\HH )$ is a smooth manifold.
\item The group $\KK_1$ acts from the left on this manifold,
by the following formula:
$$k\cdot[c_1,c_2, \dots ,c_m]=[kc_1, c_2, \dots ,c_m].$$
\item There is an inductive formula
$$\EE (\KK_1 ,\dots ,\KK_m ;\HH ) \cong 
\KK_1 \times_{\HH} \EE (\KK_2 ,\dots , \KK_m ;\HH )$$
where the right action of $\HH$ on $\KK_1$ is by right multiplication
and the left action of $\HH$ on $\EE (\KK_2 , \dots , \KK_m ;\HH )$
is via 2. In particular, there is a fiber bundle
$\EE (\KK_1 ,\dots ,\KK_m ;\HH ) \to \KK_1 /\HH $
with fiber $\EE (\KK_2 ,\dots , \KK_m ;\HH )$.
\item The fiber bundle has a section 
$$s: \KK_1 / \HH \to \EE (\KK_1 ,\dots ,\KK_m ;\HH )\quad ; \quad
k\HH \mapsto [k,e,\dots ,e].$$ 
\end{enumerate}
\end{lemma}

\begin{proof}
An easy computation shows that the action from definition 
\ref{definition:K-cycle} is free. The group
$\HH^m$ is compact. So 1. follows by
Theorem \ref{thm:sm.quotient}.

One easily checks that the action in 2. is well-defined.
In 3. one checks that the identity on $\KK_1 \times \dots \times \KK_m$
induces maps in both directions. Finally, 
$(kh, e, \dots ,e) = (k,e,\dots ,e)*(h,\dots ,h)$ for $h$ in $\HH$, 
which shows that the section in 4. is well-defined.
\end{proof}

Now assume that $M$ is a compact and connected Riemannian manifold. Then
$\GG = I_0(M)$ is a compact Lie group by the Myers-Steenrod theorem. 
The isotropy group $\GG_p$ at $p\in \MM$ is a closed subgroup.

Let $\geo : [0,1] \to M$ be a geodesic loop. Choose 
$0<t_1<t_2<\dots <t_m<1$ such that 
$\geo (t_1) , \dots , \geo (t_m)$ are the points conjugate to
$\geo (0)$ along $\geo$. 

\begin{definition}
\label{def:geod.grp}
Define closed subgroups of $\GG$ as follows:
\begin{align*}
& \KK = \GG_{\geo (0)}, \quad \KK_i = \GG_{\geo (0)} \cap \GG_{\geo (t_i)} 
\text{ for } 1\leq i \leq m, \quad \\
& \HH = \{ g\in G | g\gamma (t) = \gamma (t) \text{ for all } t \} .
\end{align*}
Define a corresponding $K$-cycle by 
\[
\kgcycle \geo \GG = \EE (\GG ,\KK_\bullet ; \HH )=\EE (\GG ,\KK_1 , \KK_2,\dots,\KK_m ; \HH)
\] 
and write $\kcycle \geo = \EE (\KK_\bullet ; \HH )$ for the fiber. 
\end{definition}

We will now introduce the Bott-Samelson map. Note that
our formula is different from the one in \cite{Ziller} page 14. We have 
an extra $\GG$ factor. 

\begin{definition}
$\widetilde \Bott_\geo : \GG \times \KK_1 \times \dots \times \KK_m
\to \Lambda M$ 
is the map given by
$$
\widetilde \Bott_\geo (g,c_1,\dots ,c_m)(t) = 
\begin{cases}
g\geo (t) &, 0\leq t \leq t_1 \\
gc_1\geo (t) &, t_1 \leq t \leq t_2\\
\vdots & \vdots\\
gc_1\dots c_{m-1}\geo (t) &, t_{m-1}\leq t \leq t_m \\
gc_1\dots c_{m-1}c_m\geo (t) &, t_m \leq t \leq 1.
\end{cases}
$$ 
\end{definition} 

Note that the geodesic pieces fit together such that 
$\widetilde BS_\geo$ takes values in $\Lambda M$ as stated. 

\begin{proposition}
\label{prop:Bott-Samelson}
The map $\widetilde \Bott_\geo$ is constant on $\HH^{m+1}$-orbits
so it induces a map
$$\Bott_\geo : \kgcycle \geo \GG \to \Lambda M$$
which is called the Bott-Samelson map. This map is $\GG$-equivariant
where $\GG$ acts from the left on $\kgcycle \geo \GG$
via lemma \ref{lemma:K-cycle}, 2. and from the left on $\Lambda M$ by
$(g, \eta )\mapsto g\eta $.
\end{proposition}

\begin{proof}
This follows by an easy direct computation.
\end{proof}

The fiber $\kcycle \geo$ is mapped to piecewise geodesic loops 
with the same initial point $\geo (0)$ as $\geo$ 
by the Bott-Samelson map. The section of the fiber bundle 
$s: \GG /\HH \to \kgcycle \geo \GG $
defines a submanifold of $\kgcycle \geo \GG $. The image
of this submanifold under the Bott-Samelson map is $\GG \geo$, 
that is the critical submanifold containing $\geo$.

Any other curve in the image of the map is not a geodesic,
but a broken geodesic. In particular, the
other critical points for the energy functional 
is the image of the translates of $\geo$ itself.
Note that all closed curves in the image of $\Bott_\geo$ 
will have the same energy since segments of $\gamma$ are acted upon by 
isometries.

\subsection{A deformation of the Bott-Samelson map and the normal bundle}
\label{sec:deformation}

We can deform the map $\Bott_\geo$ a little, as done in
\cite{Ziller} \S 3, or \cite{BS} \S 10. We fix a
suitably small positive number $\epsilon$. In dependence of
this number, we  change the curves in the vicinity of a a corner
$\geo (t_i)$. Bott and Samelson replace the curve between 
$\geo (t_i-\epsilon)$ and $\geo (t_i+\epsilon)$ with
the unique shortest geodesic between them. This
exists, because we did chose $\epsilon$ small enough.

We obtain a new map
\[
\overline \Bott_\geo : \kgcycle \geo \GG \to \Lambda \MM.
\]
It agrees with $\Bott_\geo$ on the geodesics, that is on the
image of the section $s$.

Bott and Samelson do not discuss the parametrization of the curve
one obtains, since they are interested in the length functional $L$.
We use the energy functional $E$, so the parametrization matters to us. We 
parametrize all the curves proportional to arc length. Thus, a curve 
$\sigma : [0,1]\to \MM$ has $|\sigma^\prime (t) | = k$ for a constant $k$ and all $t$.
In consequence $E(\sigma )=\frac 1 2 L(\sigma )^2$.

If $x \in \kgcycle \geo \GG \setminus s(\GG /\HH)$, then 
$\overline \Bott_\geo (x)$ will be a closed
curve with strictly lower energy than $\gamma$. Moreover,
the length functional composed with $\overline \Bott_\gamma $ 
restricted to $\kcycle \geo$ takes on a non-degenerate maximum at 
$[e,\dots ,e]$, \cite{Ziller} \S 3, \cite{BS} \S 10.  
Because of the energy-length relation above, 
one also has a non-degenerate maximum for the energy functional. 

A formula on page 14 of \cite{Ziller} states that
\begin{equation}
\label{eq:dimension}
\lambda (\gamma )=\sum_{i=1}^m \dim \KK_i /\HH,  
\end{equation}
where $\lambda (\gamma )$ is Ziller's notation for the dimension of the
negative bundle. We have used our notation on the right hand side.
By this formula we have:

\begin{theorem}
\label{dimension.negative.bundle}
The dimension of $\kcycle \geo$ agrees with the dimension 
of the negative bundle.
\end{theorem}

We can now prove an addendum to Theorem \ref{th.filtration1}.

\begin{theorem}
\label{th.filtration}
Let $\MM$ be a compact symmetric space.
For two consecutive steps $F^{\nu^\prime} \subseteq F^\nu$
in the filtration of $\Lambda \MM$ and $\gamma \in \nu$,  
we have a diagram which commutes up to homotopy
\[
\xymatrix@C=1cm{
& \kgcycle \geo \GG \ar[r]^-{c} \ar[d]_-{\overline \Bott_\geo} 
& \Th (\norbd^-_\nu ) \ar[d]^-{\simeq } \\
F^{\nu^\prime} \ar[r] & F^\nu \ar[r] & F^\nu /F^{\nu^\prime }.\\
}
\]
The map $c$ is the Thom collapse map that corresponds to the
embedding of the critical manifold of geodesic loops
$\GG /\HH \hookrightarrow \kgcycle \geo \GG$. 
\end{theorem}

\begin{proof} 
We apply lemma \ref{lemma:negative.bundle} in the following setting:
Let $\eta$ denote the normal bundle for the embedding 
$s: \GG /\HH \hookrightarrow \kgcycle \geo \GG$.
Let $f:D^o (\eta ) \to \Lambda \MM$ be the composite which first sends $D^o (\eta )$ to a 
corresponding tubular neighborhood of $s(\GG / \HH )$ in $\kgcycle \geo \GG$ and then by the 
restriction of $\overline \Bott_\geo$ maps this neighborhood to $\Lambda \MM$.
Then the zero section of $\eta$ composed with $f$ is the inclusion of the critical submanifold
$s(\GG / \HH )$ in $\Lambda \MM$ as required. The dimension condition is satisfied since
\[ \dim \eta = \dim \kgcycle \geo \GG - \dim \GG / \HH = \dim \kcycle \geo = \dim \xi^- , \]
where the last equality follows from theorem \ref{dimension.negative.bundle}.

The observation regarding the non-degenerate maximum for 
$E\circ \overline \Bott_\geo |_{\kcycle \geo}$ implies the required condition for the 
Hessian of $E\circ f$. The smoothness condition also hold, so the lemma applies and 
the result follows.
\end{proof}

To get further, we want to determine the normal bundle of $\GG/ \HH$ in 
$\kgcycle \geo \GG$.
We now return to the general setting from 
Definition \ref{definition:K-cycle}.
The point $\underline e =[e,\dots ,e]\in \EE (\KK_\bullet ;\HH )$ 
is a fixed point of the left action of $\HH$ via lemma
\ref{lemma:K-cycle}, 2. In particular, the tangent space
of $\EE (\KK_\bullet ;\HH )$ at this point is an $\HH$ representation.

\begin{theorem}
\label{th.tangent.space}
Let $\mathfrak{h}$, $\mathfrak{k}_i$ denote the Lie algebras of
the groups $\HH$, $\KK_i$ for $i=1, 2, \dots , m$. Consider these as
$\HH$-representations with the adjoint action. 
Then the $\HH$-representation 
$V:=T_{\underline e} \EE (\KK_\bullet ; \HH )$
is equivalent to the $\HH$-representation
\[
\mathfrak{k}_1 /\mathfrak{h} \oplus
\mathfrak{k}_2 /\mathfrak{h} \oplus 
\dots \oplus \mathfrak{k}_m /\mathfrak{h} .
\]
Furthermore, the normal bundle of 
$s(\GG /\HH ) \subseteq \EE(\GG ,\KK_\bullet ;\HH )$ is the
$\GG$-vector bundle $\GG \times_{\HH } V \to \GG /\HH$.
\end{theorem}

\begin{proof} 
We prove the first statement by induction on $m$.
If $m=1$, we ask for the tangent space at $[e]$
of $\KK/\HH$. 

Let $\tilde\KK$ be $\KK$ considered as a manifold with the conjugation 
left action of $\HH$. The quotient map
$\pi : \tilde \KK\to \KK/\HH$ given by $k\mapsto k\HH$ is $\HH$-equivariant and
by theorem \ref{thm:sm.quotient} it is a submersion. Consider its differential at
the unit $e\in \tilde \KK$. Since $e$ is a fixed point, the source of the differential is a 
$\HH$-representation. By definition, it is exactly $\mathfrak{k}$, with the adjoint representation.

Let $\tilde \HH$ be $\HH$ equipped with the conjugate left action of $\HH$.
The inclusion map $\iota: \tilde \HH \hookrightarrow \tilde \KK$ is $\HH$ equivariant and
it is an immersion since $\HH$ is a closed subgroup of $\KK$ by Cartan's theorem. 
The composite $\pi \circ \iota$ is trivial. Taking differentials and considering 
the dimensions involved, we obtain a short exact sequence of $\HH$-representations
$0 \to \mathfrak{h} \to \mathfrak{k} \to T_e(\HH / \KK) \to 0$. This finishes the proof 
of the induction start.

The induction step is a consequence of the following general 
argument: Let $\KK$ be a compact Lie group with a closed subgroup
$\HH$. Suppose that $\HH$ act smoothly from 
the left on a manifold $\MM$. 
Assume that $p\in \MM$ is an $\HH$-fixed point. 
The orbit space $\KK \times_\HH \MM$ is a left $\HH$ space and
$[e,p]$ is a fixed point for this action.
By theorem \ref{thm:tb.borel} the tangent $\HH$-representation at 
$[e,p]$ is isomorphic to $T_p(\MM ) \oplus \mathfrak{k} /\mathfrak{h}$.

This finishes the proof of the first statement of the theorem.
The second statement follows from theorem \ref{thm:slices}.
\end{proof}

\section{Splitting fixed points}
\label{sec:fixpt}

For a space with a group action, a special type of fixed points will produce 
a splitting of the space after a number of suspensions. We call these fixed points 
{\em splitting fixed points}. In this section we will use them to obtain 
a splitting result for free loop spaces of symmetric spaces.

\subsection{Definition and relation to Thom spaces}

Consider the following general situation.
Let $\HH$ be a compact Lie group which acts from the left on a 
based space $X$ (such that the base point is a fixed point).
Assume that $x\in X$ is a fixed point different from the basepoint. 
Also assume that there is an $\HH$-invariant open neighborhood 
$U_x$ of $x$ such that $U_x$ is  diffeomorphic to a manifold 
and the restriction of the action map $H\times U_x \to U_x$ 
is smooth. 

If $X=\MM_+$ (that is $\MM$ with a disjoint basepoint adjoined) where 
$\MM$ is a smooth manifold with a differentiable action of $\HH$
and a fixed point $x\in \MM$, then the assumptions are satisfied by 
Theorem \ref{thm:slices}, 1. This is our main example. But we want 
to keep the option of changing the topology of $\MM$ away from the fixed point.

The tangent space $T_xX$ makes sense, and is a representation  
of $\HH$. Again by Theorem \ref{thm:slices}, 1. 
there is a neighborhood $U_x$ of $x$ which is 
equivariantly diffeomorphic to $T_x X$. Collapsing the complement of 
this neighborhood defines an $\HH$-equivariant map from $X$ to
the one point compactification of the tangent space:
\[
c_x: X \to (T_x X)^+.
\]
This collapse map is independent of the choice of $U_x$ up to homotopy.

\begin{definition} 
We say that $x$ is a splitting fixed point for $X$ up to $m$-fold 
suspension, if the suspended map $\Sigma^m c_x$ is a surjective 
retraction in the $\HH$-equivariant category. That is, if there 
is a map $s: \Sigma^m (T_x X)^+ \to \Sigma^m X$ such that
$(\Sigma^m c_x) \circ s$ is equivariantly homotopic to the 
identity. 

We say that $x$ is a stably splitting fixed point for $X$ if 
there is a finite dimensional $\HH$-representation 
$V$ such that the equivariant map $V^+ \wedge c_x$ is a 
surjective retraction in the $\HH$-equivariant category. 
\end{definition}

The main example that we will consider is 
$\EE (\KK_\bullet ;\HH )$ with the left $\HH$ action
specified by lemma \ref{lemma:K-cycle}, 2. and where 
$x=\underline e = [e,\dots ,e]$. We are going to give conditions 
on the groups $\KK_i$ and $\HH$ that ensure that $\underline e$
is a splitting fixed point stably or up to suspension.

Assume that $\HH$ acts smoothly from the left on the manifold
$\MM$, and that $\HH$ is a closed subgroup of another compact Lie group $\GG$.
Then, $\GG \times_\HH \MM$ is a left $\GG$ space with the action 
$g^\prime \cdot [g,m]=[g^\prime g, m]$. If $x\in \MM$ is a fixed point for
the action of $\HH$, Theorem \ref{thm:slices}, 2. applies: We have an embedding
$\GG /\HH \hookrightarrow \GG \times_\HH \MM$ defined by $g\HH \mapsto [g,x]$ and 
the associated normal bundle $\nu (\GG /\HH )$ is the same as the $\GG$-vector bundle
$\GG \times_{\HH} (T_x \MM) \to \GG /\HH$.

\begin{lemma}
\label{lem:thsplit}
Suppose that $x \in \MM$ is a splitting fixed point for $\MM_+$ up to $m$-fold 
suspension. Then, the $m$-fold suspension $\Sigma^m c$ of the Thom collapse map
$$c: (\GG \times_\HH \MM )_+ \to \Th (\nu (\GG /\HH ))$$
is split up to homotopy.
\end{lemma}

\begin{proof} 
First note that for any pointed $\HH$-space $X$ and finite dimensional $\GG$-representation
$W$ there is a homeomorphism

\begin{equation}
\label{equiv-homeo}
\xymatrix@C=1.0 cm{
\GG_+ \wedge_\HH (W^+ \wedge X) \ar[r]^-{\cong} & W^+ \wedge( \GG_+ \wedge_\HH X) . }
\end{equation}
It is defined by $[g,w,x] \mapsto [gw,g,x]$ with inverse $[v,g,x] \mapsto [g, g^{-1}v, x]$.
Both tranformation rules are seen to be well-defined.

The Thom collapse map $c$ can be identified with the map
\[
\xymatrix@C=1.0cm{
\GG_+ \wedge_\HH c_x : \GG_+ \wedge_\HH \MM_+ \ar[r] & \GG_+ \wedge_\HH (T_x\MM )^+ .}
\]
Since $x$ is a splitting fixed point op to $m$-fold suspension, we have a splitting
map $s$ for $\Sigma^m c_x$. Via the homeomorphism (\ref{equiv-homeo}), for $W$ 
a trivial $m$-dimensional $\GG$-representation, we get a map
\begin{align*}
& \Sigma^m (\GG_+ \wedge_\HH (T_x\MM )^+ )
\cong \GG_+ \wedge_\HH \Sigma^m (T_x\MM )^+ 
\xrightarrow{\GG_+ \wedge_\HH  s} \\
& \GG_+ \wedge_\HH \Sigma^m (\MM_+) 
\cong \Sigma^m (\GG_+ \wedge_\HH \MM_+)
\end{align*}
which splits $\Sigma^m c$. 
\end{proof}

There is a corresponding statement in the stable case. 
But to prove it, we need 
the following neat result about representations of 
compact Lie groups.

\begin{theorem}
\label{thm:duister}
Let $\GG$ be a compact Lie group with a closed subgroup $\HH$.
Suppose that $V$ is a finite dimensional representation of $\HH$.
There is a finite dimensional representation $W$ of $\GG$, such that
$V$ is a subrepresentation of the restriction of $W$ to $\HH$.
\end{theorem}

For a proof of this, see for instance \cite{DK}, corollary 4.7.2.
The main ingredient of this proof is the Peter-Weyl theorem.
Since representations of compact Lie groups are semisimple,
we can actually find a finite dimensional  $\HH$-representation 
$U$ such that $U\oplus V$ is the restriction of a $\GG$-representation.

Here is the stable splitting result.
\begin{lemma}
\label{lem:stthsplit}
Assume that $x\in \MM$ is a stably splitting fixed point for $\MM_+$.
Then there exists a natural number $m$ such that the $m$-fold suspension $\Sigma^m c$
of the Thom collapse map
\[
c : (\GG \times_{\HH} \MM)_+ \to \Th (\nu (\GG /\HH ))
\]
is split up to homotopy.
\end{lemma}

\begin{proof} 
Let $s: V^+ \wedge (T_x \MM )^+ \to V^+ \wedge \MM_+$ be a 
splitting map for $V^+ \wedge c_x$. Pick a finite dimensional $\HH$-representation $U$ 
(using theorem \ref{thm:duister}), such that $U\oplus V$ is 
the restriction of a $\GG$-representation. Then 
\[
U^+ \wedge s: U^+ \wedge V^+ \wedge (T_x \MM )^+ \to 
U^+ \wedge V^+ \wedge \MM_+
\]
is a splitting map for $U^+ \wedge V^+ \wedge c_x$ which gives a map
\begin{equation}
\label{eq:splitting.map}
\GG_+ \wedge_\HH (U\oplus V)^+ \wedge (T_x\MM )^+
\to 
\GG_+ \wedge_\HH (U\oplus V)^+ \wedge \MM_+ 
\end{equation}
Non-equivariantly we have $S^m\cong (W\oplus V)^+$ for some $m$. 
Applying the homeomorphism (\ref{equiv-homeo}) to the source and to the 
target of the map (\ref{eq:splitting.map}) gives us the desired splitting
\[
\Sigma^m \Th (\nu (\GG /\HH ) ) \cong 
\Sigma^m \GG_+ \wedge_\HH (T_x\MM )^+ \to 
\Sigma^m \GG_+ \wedge_\HH \MM_+ \cong
\Sigma^m (\GG \times_\HH \MM)_+.
\] 
\end{proof}

\subsection{Existence of splittings}

Here is an elementary example of a fixed point which provides a splitting after a single
suspension.

\begin{example}
\label{ex.repr} 
Let $V$ be a real vector space with $2\leq \dim V < \infty$ and equipped with an 
inner product. Let $\HH$ be a closed subgroup of the orthogonal group $O(V)$. The
canonical action of $O(V)$ on the unit sphere $S(V)$ in $V$ restricts to 
an action $\HH \times S(V) \to S(V)$. Assume that $x\in S(V)$ is a fixed point for this action. 
Then $x$ is a splitting fixed point for $S(V)_+$ up to $1$-fold suspension.

The argument is as follows:
The subspace $L=\RR x \subseteq V$ has trivial $\HH$-action.
Its orthogonal complement $U=L^\perp \subseteq V$ is an $\HH$-invariant subspace
(follows from the equation $\langle hu,x \rangle = \langle hu, hx \rangle = \langle u, x \rangle$
for $h\in \HH$). Consider stereographic projection
\[ \pi : S(L\oplus U) \to U^+ \]
which maps $-x$ to $\infty$. It is a homeomorphism. 
The line through $-x$ and $v$ has parametric equation
$\ell (t) = -x+t(v+x)$. Solving 
$\langle \ell (t),x \rangle = 0$ for $t$ and inserting the solution, one finds the
transformation rule
\[ \pi (v) = -x + \frac {v+x} {1+ \langle v, x \rangle } , \quad v\in S(V) \setminus \{ -x \} . \]
It follows by this rule that $\pi$ is an $\HH$-equivariant map since $x$ is a fixed point and 
$\HH$ is a subgroup of $O(V)$.

Thus we can use the collapse map $c_x : S(V)_+ \to (T_xS(V))^+$
given by $c_x (v) = \pi (v)$ for $v\in S(V)$ and $c_x (+) = \infty$. This map is
equivariantly homotopic to the map 
\[ {\widetilde c}_x = \pi^{-1}\circ c_x : S(V)_+\to S(V)\] 
which restricts to the identity on $S(V)$ and maps $+$ to $-x$. 
We claim that the map $\Sigma \tilde c_x$ has an equivariant section. 

To see this, we note that one can identify it with the quotient map
\[
q: S(V\oplus \epsilon )/\{ \pm (0,1)\}
\to S(V\oplus \epsilon )/\{ (-ax,b) \mid a^2+b^2=1, a\geq 0\} ,
\]
where $\epsilon$ denotes a trivial one dimensional 
$\HH$-representation. Consider the following quotient map:
\[
p: S(V\oplus \epsilon ) \to 
S(V\oplus \epsilon )/\{ \pm (0,1) \}.
\]
The composite $q\circ p$ of the two maps is a quotient map given by
dividing out a contractible subspace (actually, an interval). 
Since the inclusion of the subspace is an equivariant cofibration,
this composite is an equivariant homotopy equivalence.
Composing the homotopy inverse of this homotopy equivalence
with the map $p$, we get a homotopy left inverse of
$q$ and thus also of $\Sigma\tilde c_x$ as required. 
\end{example}

Here is a splitting result of a more general nature.

\begin{lemma}
\label{lem:rep.split}
Let $\GG$ be a compact Lie group with a closed subgroup $\HH$ and
let $W$ be a finite dimensional real representation of $\HH$.
\begin{enumerate}
\item Assume that $[e]\in (\GG / \HH )_+$ is stably splitting.
If the class $[W]$ is contained in the image of the 
restriction map $RO(\GG )\to RO(\HH )$ of real representation rings,
then the $\HH$-fixed point $[(e,0)]\in \GG_+ \wedge_\HH W^+$ 
is stably splitting.
\item Assume that $[e]\in (\GG /\HH )_+$ is splitting up to $m$-fold suspension.
If the $\HH$-representation $T_{[e]}(\GG / \HH )\oplus W\oplus \epsilon^m$ 
extends to a $\GG$-representation, then
$[(e,0)]\in \GG_+ \wedge_\HH W^+$ is splitting up to $m$-fold suspension.
\end{enumerate}
\end{lemma}

\begin{proof}
1. The condition on $W$ says that there are $\GG$-representations 
$U,U^\prime$ so that after restricting the action to the subgroup
$\HH$ we have an isomorphism of $\HH$-representations 
$W\oplus U \cong U^\prime$.

The point $[e]\in (\GG / \HH)_+$
is stably splitting so there exist an $\HH$-representation $V$ 
and a splitting $s$ of the map 
$V^+\wedge c_{[e]} : V^+\wedge (\GG / \HH)_+ \to 
V^+\wedge (T_{[e]}(\GG / \HH))^+ .$
For any manifold $\MM$ and vector space $E$ one has an
isomorphism for $x\in \MM$ as follows:
\begin{equation}
\label{tangentid}
T_{(x,0)} (M_+\wedge E^+) \cong T_{(x,0)} (M\times E) \cong T_x(M) \oplus E.
\end{equation}
We get a commutative diagram of $\HH$-equivariant maps
\[
\xymatrix@C=1cm{
(\GG /\HH )_+ \wedge W^+ \wedge U^+ \ar[r]^-{c_{([e],0,0)}} 
\ar[rdd]_-{c_{[e]} \wedge W^+ \wedge U^+} &
\big( T_{([e],0,0)} \big( (\GG /\HH )_+ \wedge W^+ \wedge U^+ \big) \big)^+ \ar[d]^-{\cong } \\
& \big( T_{[e]}(\GG / \HH ) \oplus W \oplus U \big)^+ \ar[d]^-{\cong } \\
& \big( T_{[e]} (\GG / \HH ) \big)^+ \wedge W^+ \wedge U^+ .
}
\]
Thus, the map $s\wedge W^+\wedge U^+$ gives an $\HH$-equivariant section of
${V^+\wedge c_{([e],0,0)}}$. 

Consider the composite of the following two $\GG$-homeomorphisms:
\[ 
\xymatrix@C=.7cm{
\theta : (\GG_+ \wedge_\HH W^+)\wedge U^+ \ar[r]_-{\cong}^-{\phi}
& \GG_+ \wedge_\HH (W^+ \wedge U^+) \ar[r]_-{\cong}^-{\psi}
& (\GG /\HH )_+ \wedge W^+\wedge U^+.
}
\]
Here $\phi ([g,w],u) = [g,(w,g^{-1}u)]$ and 
$\psi ([g,(w,u)]) = ([g], g\cdot (w,u))$. The map
$\psi$ exists since $W\oplus U$ is actually the restriction of 
a $\GG$-representation. Note that we may consider $\theta$ 
a local diffeomorphism near the point $p=([e,0],0)$ and that $\theta (p) = ([e],0,0)$. 

We have a homotopy commutative diagram, where horizontal maps
are homeomorphisms:
\[
\xymatrix@C=1.5cm{
V^+\wedge (\GG_+ \wedge_\HH W^+)\wedge U^+ 
\ar[r]^-{V^+\wedge \theta} \ar[d]^-{V^+\wedge c_p }
& V^+\wedge (\GG /\HH )_+ \wedge W^+\wedge U^+ \ar[d]^-{V^+\wedge c_{\theta (p)}} \\
V^+\wedge T_p ( (G_+ \wedge_\HH W^+)\wedge U^+ )^+
\ar[r]^-{V^+\wedge (d_p \theta )^+ } 
& V^+\wedge T_{\theta (p)} ( (\GG /\HH )_+ \wedge W^+\wedge U^+ )^+ .
}
\]
The right vertical map has a splitting $s\wedge W^+ \wedge U^+$
as remarked above. So the left vertical map also has a splitting.
The result now follows by using (\ref{tangentid}) on the lower left corner of the 
diagram.

2. The $m$-fold suspension case is proved in the same way, 
using the trivial representation $\epsilon^m$ in place of $U$.
\end{proof}
 
We can now prove the main result on existence of splittings.

\begin{theorem}
\label{thm:susp.split}
Let $\GG$ be a compact Lie group, and $\HH\subseteq \GG$ a closed
subgroup. Assume that $\HH$ act smoothly on a manifold $\MM$ with a fixed point $x$.
\begin{enumerate}
\item If $[e] \in (\GG /\HH )_+$ and $x \in \MM_+$ are both stably 
splitting $\HH$-fixed points and the class $[T_xM]$ is in the image of the restriction
$RO(\GG )\to RO(\HH )$ then
$[e,x]\in (\GG \times_\HH \MM)_+$ is also a stably splitting $\HH$-fixed point.
\item If $[e]\in (\GG / \HH)_+$ and $x\in M_+$ are 
both splitting $\HH$-fixed points up to $m$-fold suspension, and if the $\HH$-representation
$T_xM\oplus \epsilon^m$ extends to a $\GG$-representation, 
then $[e,x]\in (\GG \times_\HH \MM)_+$ is also a splitting $\HH$-fixed point up 
to $m$-fold suspension.
\end{enumerate}
\end{theorem}

\begin{proof}
1. Let $W=T_xM$ be the tangent space representation of $\HH$. 
The collapse map $\MM_+\to W^+$ is a diffeomorphism
in a neighborhood of $x$, so we have a diagram
\[
\xymatrix@C=1cm{
\GG_+\wedge_\HH \MM_+  \ar[r]^-{\GG_+ \wedge_\HH c_x} \ar[d]^-{c_{[e,x]}}
& \GG_+\wedge_\HH W^+ \ar[d]^-{c_{[e,0]}} \\
(T_{[e,x]}(\GG \times_\HH \MM ))^+ \ar[r]^-{\cong } 
& (T_{[e,0]}(\GG_+ \wedge_\HH W^+))^+ .  
}
\] 
The right vertical map is stably split by lemma \ref{lem:rep.split}.
The upper horizontal map is stably split, since $c_x$ is by assumption.
It follows that the left vertical map is stably split. 

2. If $c_x$ is $m$-fold suspension split, we see that 
the right and upper maps are split in the diagram.
\[
\xymatrix@C=1cm{
S^m\wedge \GG_+ \wedge_\HH \MM_+ \ar[r] \ar[d]  
& S^m\wedge \GG^+ \wedge_\HH W^+ \ar[d] \\
S^m\wedge (T_{[e,x]} (\GG \times_\HH \MM))^+ \ar[r]^-{\cong } 
&S^m\wedge (T_{[e,0]} (\GG_+ \wedge_\HH W^+))^+ \\  
}
\]
The result follows.
\end{proof}

\subsection{Relations to stable equivariant framings}

A slightly different approach to the existence of stable
splittings is to construct them from stable framings.
In the non-equivariant case, it is known that the existence of 
a stable splitting of the collapse map 
$\MM_+ \to (T_x\MM )^+$ is equivalent to the existence
of a stable fiber homotopy trivialization of the  
tangent sphere bundle.  We use an equivariant version
of this idea.

Assume that the compact Lie group $\HH$ acts smoothly
on the compact manifold $\MM$. 
If $V$ is a representation of $\HH$ and $X$ is an $\HH$-space,
let $\epsilon^V (X)$ denote the trivial $\HH$-vector bundle 
with projection map $\mathop{pr}_X : X\times V\to X$. 

\begin{definition}
We say that $\MM$ is stably framed if there are $\HH$-representations 
$V$ and $W$, such that there is an 
equivalence of $\HH$-vector bundles 
\[ \tau (\MM ) \oplus \epsilon^V (\MM ) \cong \epsilon^W (\MM ). \]
\end{definition}

\begin{lemma}
\label{lem:framed.split}
If $\MM$ is stably framed, and $x\in \MM$ is any fixed point,
then $x\in M_+$ is a stably splittting fixed point.
\end{lemma}

\begin{proof}
We first show that there exists an equivariant embedding $\MM \hookrightarrow U$ into some
finite dimensional $\HH$-representation $U$ such that the associated 
normal bundle $\nu (\MM )$ is a trivial $\HH$-vector bundle i.e. 
$\nu (\MM ) \cong \epsilon^V (\MM )$ for a finite dimensional $\HH$-representation $V$.

There is an equivariant embedding $\MM \hookrightarrow \overline U$ 
into an orthogonal finite dimensional $\HH$-representation as in \cite{Bredon} theorem VI.4.1. 
Write $\overline \nu (\MM )$ for the associated normal bundle. We have 
$\tau (\MM )\oplus \overline \nu (\MM ) \cong \epsilon^{\overline U}(\MM )$.

Since $\MM$ is stably framed, there exist $\HH$-representations $\overline V$ and $\overline W$
such that $\tau (\MM ) \oplus \epsilon^{\overline V} (\MM ) \cong \epsilon^{\overline W} (\MM )$.
Let $U=\overline U \oplus \overline W$ and consider the embedding
$\MM \hookrightarrow \overline U \hookrightarrow U$ with associated normal bundle $\nu (\MM )$.
We have
\[ \nu (\MM ) \cong \overline \nu (\MM ) \oplus \epsilon^{\overline W}(\MM ) \cong 
\overline \nu (\MM ) \oplus \tau (\MM ) \oplus \epsilon^{\overline V} (\MM )  \cong
\epsilon^{\overline U}(\MM ) \oplus \epsilon^{\overline V}(\MM ) . \]
Thus, $\nu (\MM ) \cong \epsilon^V (\MM )$ where $V=\overline U  \oplus \overline V$.
So we have the desired embedding.

The Thom collapse map 
$U^+ \to \Th (\nu (\MM ) )$ gives us an $\HH$-equivarinat map $s$ as follows:
\[
V^+ \wedge (T_x\MM )^+ \cong (V\oplus T_x\MM )^+ \cong U^+ \to 
\Th (\nu (\MM ) )\cong \Th (\epsilon^V (\MM )) \cong V^+\wedge \MM_+.
\]
The composite $(V^+ \wedge c_x) \circ s$ maps $[v, w]$ to $[v,0]$ where $v\in V$ and $w\in T_x \MM$.
So it is equivariantly homotopic to the identity through the homotopy $H([v,w],t)=[v,tw]$, $0\leq t\leq 1$. 
\end{proof}

Suppose that $\HH$ is a closed subgroup of another compact
Lie group $\KK$. The left $\KK$-action on $\KK /\HH$ restricts to a
left $\HH$-action.

\begin{lemma}
\label{lemma:representation.ring}
Assume that $\KK /\HH$ is stably framed as an $\HH$-space. Also suppose that $\MM$ is stably 
framed as an $\HH$-space, such that we have an equivalence 
$\tau (\MM ) \oplus \epsilon^V (\MM ) \cong \epsilon^W (\MM )$.
Assume finally that the element $[W]-[V]\in RO(\HH )$ is contained
in the image of the restriction map $RO(\KK ) \to RO(\HH )$. Then $\KK \times_\HH \MM$ is  
stably framed as an $\HH$-space, where the action of $\HH$ is given as the restriction of
the $\KK$-action.
\end{lemma}

\begin{proof}
By theorem \ref{thm:tb.borel} there is an isomorphism of $\KK$-vector bundles
\[ \tau (\KK \times_\HH M) \cong \pi^* (\tau (\KK / \HH )) \oplus (\KK \times_\HH \tau (M)) \]
where $\pi : \KK \times_\HH \MM \to \KK / \HH$ is the projection map $\pi ([k,x]) = k\HH$.
We consider the two summands separably. 

By assumption, there exist $\HH$-representations $A$ and $B$ and an 
equivalence $\tau (\KK /\HH ) \oplus \epsilon^A (\KK /\HH ) \cong \epsilon^B (\KK /\HH )$.
We pull back this equivalence along the equivariant map $\pi$ and get the desired
result for the first summand.

From the equivalence $\tau (\MM ) \oplus \epsilon^V (\MM ) \cong \epsilon^W (\MM )$ we
get an equivalence of $\KK$-vector bundles 
\[  (\KK \times_\HH \tau (\MM ))\oplus (\KK \times_\HH \epsilon^V (\MM )) \cong
\KK \times_\HH \epsilon^W (\MM ) . \]
However, for any $\KK$-representation $U$, the $\KK$-vector bundle 
$\KK \times_\HH \epsilon^U (\MM)$ is
trivial. A trivialization is given by
\[  \KK \times_\HH (M\times U) \xrightarrow{\cong } (\KK \times_\HH M)\times U ;
\quad [k, (m,u)] \mapsto ([k,m],ku). \]
So it suffices to see that we can replace $V$ and $W$ by representations which are restrictions
of $\KK$-representations.

By Theorem \ref{thm:duister}, there exists an $\HH$-representation $V^\prime$
such that $V\oplus V^\prime$ is the restriction of an $\KK$-representation.
After replacing $V$ by $V\oplus V^\prime$ and $W$ by $W\oplus V^\prime$ we
can assume that $V$ is the restriction of a $\KK$-representation. By assumption,
$[W]-[V]=[U^\prime]-[U]$ where $U$ and $U^\prime$ are restrictions of $\KK$-representations.
Thus $W\oplus U\cong V\oplus U^\prime$ is the restriction of a $\KK$-representation and
we have the desired equivalence 
$\tau (\MM ) \oplus \epsilon^{V\oplus U}(\MM ) \cong \epsilon^{W\oplus U} (\MM )$.
\end{proof}

\begin{remark} 
The condition on the representations $V$, $W$ in the lemma is needed
as the following example shows: Let $\HH =C_2=\langle T \rangle$ be the cyclic group
of order two inside the unit circle group $\KK =S^1$.  Let $W$ be the orthogonal $\HH$-representation
given by the reflection action $T\cdot (x,y)=(x,-y)$ on $\RR^2$, where $\RR^2$ is equipped with 
the standard inner product. Let $\MM \subseteq W$ be the unit circle. 
The tangent bundle $\tau (\MM )$ consists of
pairs $(p,v)\in \MM \times W$ such that $p$ and $v$ are orthogonal vectors. 
The normal bundle $\nu (\MM )$ consist of pairs $(p,w) \in \MM \times W$ such
that $w=cp$ for some real scalar $c$. We have an equivalence of $\HH$-vector bundles
$\tau (\MM )\oplus \nu (\MM ) \cong \epsilon^W (\MM )$. But there is a trivialization
of the normal bundle $\epsilon^V (\MM) \cong \nu (\MM)$; $(p,s)\mapsto (p,sp)$ 
where $V=\RR$ with trivial $\HH$-action. Thus, $\MM$ is stably 
$\HH$-framed. 
However, if we we view $\KK$ and $\MM$ as unit circles in the complex plane, the $\HH$-actions
are given by $T\cdot e^{it} = e^{i(t+\pi )}$ and $T\cdot e^{it}=e^{i(-t)}$ respectively.
So the quotient space $\KK \times_\HH \MM$ is a Klein bottle and hence not orientable. 
This means that it cannot be even non-equivariantly stably framed.
\end{remark}

\begin{remark}
Homogeneous spaces $\KK /\HH$ are not always stably framed.
For instance, most projective spaces like $\RP^n$, $\CP^n$ or $\HP^n$  
have non trivial Stiefel Whitney classes, 
cf. \cite{Characteristic}, especially Corollary 11.15.
So they cannot even be non-equivariantly stably framed.
\end{remark}

\subsection{Application to the Bott-Samelson map}
\begin{theorem}
\label{thm:thom.split}
Let $\HH \subseteq \KK_i \subseteq \GG$ be compact Lie groups for
$1\leq i\leq m$. The normal bundle $\nu (\GG /\HH )$ of the submanifold
$s(\GG /\HH )$ in $\EE (\GG ,\KK_\bullet ;\HH )$ is the same as the 
$\GG$-vector bundle $\GG \times_\HH V \to \GG /\HH$, where
$V = \oplus_{s=1}^{m} (\mathfrak{k}_s / \mathfrak{h}) $.
We have an associated Thom collapse map 
\[ c: \EE (\GG ,\KK_\bullet ;\HH )_+ \to \Th (\nu (\GG /\HH ) ). \]
Define $\HH$-representations $V_j$ for $2\leq j \leq m$ by
$V_j=\oplus_{s=j}^{m} (\mathfrak{k}_s / \mathfrak{h} )$.  
\begin{enumerate}
\item Assume that $\KK_i /\HH$ is stably framed as an $\HH$-manifold for each $i$. 
Also suppose, that the class $[V_j] \in RO(\HH )$ is in the 
image of the restriction map $RO(\KK_{j-1} ) \to RO(\HH )$ for $2\leq j \leq m$.
Then, there is some integer $n$ such that the
$n$-fold suspended collapse map $\Sigma^n c$ splits.
\item Assume that there is an integer $k$
such that $[e]\in (\KK_i /\HH)_+$ is splitting up to $k$ fold suspension for each $i$, and that
$V_j \oplus \epsilon^k$ extends to a representation of $\KK_{j-1}$
for $2\leq j \leq m$. Then, the $k$-fold suspended collapse map $\Sigma^k c$ splits.
\end{enumerate}
\end{theorem}

\begin{proof}
The statement about the normal bundle is part of the conclusion of theorem \ref{th.tangent.space}.

1. Consider $\EE [i] =\EE (\KK_i ,\KK_{i+1} ,\dots ,\KK_m ;\HH )$.  
We claim that $\EE [i]$ is stably framed as an $\HH$-space for 
$1\leq i \leq m$. 
The argument is by downwards induction. For 
$i=m$ we have that $\EE [m] =\KK_m /\HH $, which is stably framed by 
assumption. Assume inductively that $\EE [j]$ is stably framed. 
This means that there are representations 
$W$ and $U$ of $\HH$, such that we have an isomorphism of 
$\HH$-vector bundles
\[
\tau (\EE [j]) \oplus \epsilon^W (\EE [j]) \cong 
\epsilon^U (\EE [j]).
\]
Since $\underline e\in \EE [j]$ is a fixed
point, the restriction of this isomorphism to the fibers over $\underline e$ gives us an   
isomorphism of $\HH$-representations 
$T_{\underline e} \EE [j]\oplus W\cong U$.
But according to theorem \ref{th.tangent.space} the representation 
$T_{\underline e} \EE [j]$ is exactly $V_j$. So we obtain the relation
\[
[U]-[W]=[V_j ]\in RO(\HH ).
\]
Now by our other assumption, $[U]-[V]$ is in the image of the restriction map 
$RO(\KK_{j-1})\to RO(\HH )$.
Since $\KK_{j-1} /\HH$ is assumed to be stably framed as an 
$\HH$-manifold, We can apply lemma \ref{lemma:representation.ring} 
to this situation, and get that 
$\EE [j-1] =\KK_{j-1} \times_\HH \EE [j]$ 
is indeed stably framed as an $\HH$-manifold. 

We now know that $\EE [1]$ is stably framed, so 
$\underline e \in \EE [1]$ is stably split by 
lemma \ref{lem:framed.split}. The result follows from 
lemma \ref{lem:stthsplit}.
   
2. The unstable case is treated in essentially the same way. 
We use the second part of theorem \ref{thm:susp.split} to show 
inductively that $\underline e\in \EE [j]$ is $m$-fold suspension split. 
Then the claim follows from lemma \ref{lem:thsplit}.
\end{proof}

Before we apply this result to the filtration of the free loop space, 
we show a technical lemma, which we need in the 1-fold suspension 
split case.

\begin{lemma}
\label{lem:homotopy}
Let $B$ be a connected CW complex with a connected sub complex $A$ 
and let $i: A\to B$ denote the inclusion map. Let $q$ be the
quotient map given by 
$A_+ \xrightarrow{i_+} B_+ \xrightarrow{q} B/A$.
Suppose that $\Sigma q :\Sigma (B_+)\to \Sigma (B/A)$ is a surjective
retraction up to homotopy with a right inverse up to homotopy
$\Sigma (B/A) \to \Sigma (B_+)$. Then $\Sigma (B_+)$ is
homotopy equivalent to $\Sigma(B/A)\vee \Sigma(A_+)$.
\end{lemma}

\begin{proof} 
Chose a 0-cell $a_0 \in A \subset B$ as base point.
Let $f_A :A_+\to A$ be the based map which is the identity 
outside the disjoint basepoint $+$. Let 
$\bar q : B \to B/A$ denote the quotient map.
There is a commutative diagram
\[
\begin{CD}
\Sigma (A_+) @>{\Sigma i_+}>> \Sigma (B_+) 
@>{\Sigma q}>> \Sigma (B/A)\\
@V{\Sigma f_A}VV  @V{\Sigma f_B}VV  @| \\
\Sigma A @>{\Sigma i}>> \Sigma B 
@>{\Sigma \bar q}>> \Sigma (B/A)\\
\end{CD}
\]
The bottom line is a cofibration sequence. It splits by the 
composite map
$(\Sigma f_B) \circ s : \Sigma B/A \to \Sigma B$.
In particular there is an isomorphism 
\[
H_*(\Sigma i \vee (\Sigma f_B) \circ s)\colon
H_*(\Sigma A \vee \Sigma B/A) \to H_*(\Sigma B).
\]
Since $\Sigma i\vee (\Sigma f_B)\circ s$ is a homology
equivalence between simply connected spaces,
it is a homotopy equivalence. 
So $\Sigma B\simeq \Sigma A\vee \Sigma B/A$, and obviously 
$\Sigma B \vee S^1 \simeq \Sigma A\vee S^1\vee\Sigma B/A$.

We finish the proof of the lemma by noting that 
for every connected, based CW complex $(X,x_0)$ there are two
homotopy equivalences
$$\Sigma (X) \vee S^1 \xleftarrow{\simeq } S(X) \cup_{S^0} D^1
\xrightarrow{\simeq } \Sigma (X_+)$$
which are the quotient maps for the two contractible subspaces
$I\times {x_0}$ and $D^1$ respectively. (Here $S(X)$ denotes
the unreduced suspension of $X$). 
\end{proof}

Now we return to closed geodesics on a symmetric space. For each
closed geodesic $\gamma$ , we consider the
isotropy group of the geodesic $\HH (\gamma )$. For each
point $\gamma(t_i)$, $1\leq i\leq m(\gamma )$ conjugated to 
$\gamma(0)$ along $\gamma$ we consider the group $\KK_i (\gamma )$
of isometries fixing both $\gamma (0)$ and $\gamma (t_i)$. We define
$\HH$-representations $V_j (\gamma)$ for $2\leq j \leq m(\gamma )$ as
$V_j(\gamma ) = 
\bigoplus_{s=j}^{m(\gamma )} ( \mathfrak{k}_s(\gamma )/ \mathfrak{h}(\gamma ))$.

\begin{theorem}
\label{thm:splitting}
Let $\MM$ be a connected compact symmetric space. 
Assume that all non-trivial critical submanifolds 
have positive dimensional negative bundles, that is, the manifold of constant 
paths is the set of all critical points of index zero. 
\begin{enumerate}
\item If for each closed geodesics $\gamma$, the quotient 
$\KK_i (\gamma ) /\HH (\gamma )$ is stably framed as an $\HH (\gamma )$-manifold 
for $1\leq i \leq m(\gamma )$ and the class $[V_j(\gamma )]$ is in the image of the restriction map 
$RO(\KK_{j-1} (\gamma ))\to RO(\HH (\gamma ))$ for $2\leq j \leq m(\gamma )$, then there
is a splitting of suspension spectra up to homotopy
\[
\Sigma^\infty (\Lambda M) _+ \simeq 
\Sigma^\infty M_+ \vee \bigvee_\nu \Sigma^\infty \Th (\norbd_\nu^- ).
\]
\item If there is a positive integer $k$ such that for each closed geodesic $\gamma$ one  
has that $[e] \in (\KK_i (\gamma ) / \HH (\gamma ))_+$ is splitting up to $k$-fold
suspension and $V_j(\gamma )\oplus \epsilon^k$ extends to a representation of $\KK_{j-1}(\gamma )$
for $2\leq j \leq m(\gamma )$, then there is a splitting up to homotopy
\[
\Sigma^k (\Lambda M) _+ \simeq
\Sigma^k M_+ \vee \bigvee_\nu \Sigma^k \Th( \norbd_\nu^- ).
\]
\end{enumerate}
In both cases, $\nu$ runs through the $\GG$-equivalence classes of closed geodesics of positive energy,
where $\GG = I_0(\MM )$.
\end{theorem}

\begin{proof}
We prove 2. The proof of 1. is similar but easier.
Consider the filtration of theorem \ref{th.filtration1}. According to 
theorem \ref{th.filtration}, if $\nu^\prime$ and 
$\nu$ are two consecutive steps in the filtration, we have a diagram
 \[
\begin{CD}
@. \Sigma^k ({\kgcycle \gamma \GG }_+ ) @>{\Sigma^k c}>> 
\Sigma^k \Th (\norbd_\nu ^- ) \\
@. @V{\Sigma^k \overline{BS} }VV @V f V\cong V \\
\Sigma^k F^{\nu^\prime} @>>> \Sigma^k F^\nu @>{\Sigma^k q}>> 
\Sigma^k F^\nu /F^{\nu^\prime } .\\
\end{CD}
\]
By the conclusion of theorem \ref{thm:thom.split},
$\underline e \in \kgcycle \gamma \GG $ is $k$-fold suspension 
split. So we have a splitting of $\Sigma^k c$, say 
$s: \Sigma^k \Th (\norbd_\nu ^- )\to 
\Sigma^k ({\kgcycle \gamma \GG }_+)$.
But then the map $\Sigma^k q$ is a split surjection up to
homotopy, with right inverse 
$\Sigma^k \overline{BS} \circ s \circ f^{-1}$.

Each space $\Th (\norbd_\nu ^-)$ is connected, $M$ is connected, 
so by induction, using the assumption that there are no critical 
manifolds of index 0, we see that each space $F^\nu$ is connected. 
Using lemma \ref{lem:homotopy} we see inductively that we have a 
splitting up to suspension
\[
\Sigma^k (\Lambda^{e+\epsilon} M)_+ \simeq \Sigma^k M_+ \vee 
\bigvee_{E(\nu) < e+\epsilon }
\Sigma^k \Th (\norbd_\nu^- ).
\] 
Passing to the direct limit preserves
homotopy equivalences.
\end{proof}

\section{Rank one symmetric spaces}
\label{sec:Rank1}

We will now assume that $\MM$ is {\em isotropic} (definition \ref{def:isotropic}). This condition is 
quite constraining. It is equivalent to requiring that the symmetric space has rank one. According 
to \cite{Wolf}, 8.12.2 this class includes spheres, complex and quaternionic 
projective spaces and the Cayley projective plane, but no other compact, 
simply connected symmetric spaces.

Let $\GG$ be the isometry group of $M$. 
Since we are assuming isotropy, it acts 
transitively on set of pairs $(p,v)$ where
$p\in \MM$, and $v\in T_p\MM$. It follows
that $\GG$ acts transitively on the primitive geodesics. Such a geodesic has no self intersections 
(since the tangent vector field along $\gamma$ is the restriction 
of a Killing vector field on $\MM$). So the space of primitive 
geodesics can be identified as a homogeneous space 
$\GG /\HH$ where $\HH$ is the subgroup of $\GG$ that fixes 
$\gamma$ pointwise. 

In particular, all geodesics are closed, and all primitive 
closed geodesics have the same length. Let us normalize 
the metric on $\MM$, so that they have length one. We see that we 
can also think of the space of primitive closed geodesics as the 
unit sphere bundle of the tangent bundle of $\MM$.

We take a close look at the three non-sphere types of compact,
simply connected, symmetric spaces.
We will consider all representations as representations over the 
real numbers. The dimension of a representation will thus mean its real dimension.

\subsection{The complex projective spaces}

For $\MM=\CP^n$, we can explicitly compute the various isotropy groups
mentioned above. The metric on $\CP^n$ is (a rescaling of) the Fubini-Study metric. Let us review
some of its properties. 

Consider $\CC^{n+1}$ with the standard inner product 
$\langle v, w \rangle = \sum v_j \overline{w}_j$.
The real part $g(\cdot , \cdot) = \Re \langle \cdot , \cdot \rangle$ is the standard inner product on 
$\RR^{2n+2}$, if one uses the identification $\RR^{n+1} \oplus \RR^{n+1} \cong \CC^{n+1}$; 
$(a,b) \mapsto a+ib$, and the associated norms of $g$ and $\langle \cdot , \cdot \rangle$ agree. 
Furthermore, 
\[ \langle v, w \rangle = g(v,w)-ig(iv,w). \]
Let $S^{2n+1}$ denote the unit sphere in $\CC^{n+1}$ and let 
\[ \rho : S^{2n+1} \to S^{2n+1}/ U(1) = \CP^n \]
be the canonical projection. The inner product $g$ gives a Riemannian metric on $S^{2n+1}$. For 
$x\in S^{2n+1}$ one has the orthogonal complement
\[ (\CC x)^\perp = \{ v \in \CC^{n+1} | \langle x, v \rangle =0 \}
= \{ v \in \CC^{n+1} | g(x, v)=0, \quad g(ix, v)=0 \} . \]
View it as a real subspace of $T_x(S^{2n+1}) = \{ v \in \CC^{n+1} | g(x,v)=0 \}$. 
The Fubini-Study metric is determined by the property that the composite map
\[ d\rho_x : (\CC x)^\perp \hookrightarrow T_x(S^{2n+1}) \xrightarrow{D_x (\rho )} T_{\rho (x)} (\CP^n ) \]
is an $\RR$-linear isometry for all $x$. The idendity $\rho = \rho \circ (\alpha \cdot)$ for $\alpha \in U(1)$ 
induce a commutative diagram of tangent spaces. Via this diagram one finds that
\begin{equation} \label{eqn:tangent_space_id}
 d\rho_x (v) = d\rho_y(w) \Leftrightarrow \exists \lambda \in U(1): y=\lambda x \wedge w=\lambda v .
\end{equation}

Regarding the isometry group, first note that $U(n+1)$ acts by isometries on $S^{2n+1}$ as 
$(A,x)\mapsto Ax$. The action induces a well-defined action of $U(n+1)$ on $\CP^n$ as
$(A, [x]) \mapsto [Ax]$. The linear isometry $A\cdot :T_x(S^{2n+1}) \to T_{Ax}(S^{2n+1})$
restricts to a linear isometry $(\CC x)^\perp \to (\CC Ax)^\perp$. So $U(n+1)$ acts by isometries on $\CP^n$.
In fact every isometry of $\CP^n$ is induced by a unitary matrix. The action is however not effective. Its kernel, 
that is the set of unitary matrices $A$ such that $Ap=p$ for all $p$ in $\CP^n$, is seen to be the subgroup
\[ \zun_{n+1} = \{ zI_{n+1} | z \in U(1) \} . \]
Thus, the isometry group of $\CP^n$ is the projective unitary group 
\[ I(\CP^n ) = PU(n+1) = U(n+1)/ \zun_{n+1} . \]
This group is connected so $\GG = I_0(\MM ) = PU(n+1)$.

We determine the various isotropy groups.
Let $e_1, \dots , e_{n+1}$ denote the standard basis for $\CC^{n+1}$.
Consider the point $p=[e_{n+1}]\in \CP^n$. The isotropy group 
under the $U(n+1)$ action of this point is 
$U(n+1)_p = U(n) \times U(1)$. 
We obtain the isotropy group at the point
under the action of the isometry group by factoring out $\zun_{n+1}$. 
This isotropy group is  
$\KK_2 = \GG_p = (U(n)\times U(1))/\zun_{n+1}$,
where the name of the group is chosen to be compatible
with the notation of section \ref{sec:fixpt}. There is an isomorphism 
of groups $\KK_2 \cong U(n)$; $[A, z]\mapsto z^{-1}A$ with inverse $A\mapsto [A,1]$.

Consider a second point $q=[e_n]\in \CP^n$. The subgroup of 
$U(n+1)$ that preserves both $p$ and $q$ is
$U(n+1)_p \cap U(n+1)_q = U(n-1)\times U(1)\times U(1)$. 
The corresponding group of isometries is 
$\KK_1= (U(n-1)\times U(1)\times U(1))/\zun_{n+1}$.
There is a group isomorphism $\KK_1 \cong U(n-1)\times U(1)$. 

Let $\gamma$ be a unit speed geodesics passing through $p$ and $q$. Let $[A]$ be an isometry which preserves
$\gamma$ pointwise. Equivalently, $[A]$ preserves both $p$ and the velocity vector of $\gamma$ at $p$. So by (\ref{eqn:tangent_space_id}) there exists a $\lambda \in U(1)$ such that $Ae_{n+1} = \lambda e_{n+1}$ and
$Ae_n = \lambda e_n$. Thus the subgroup that fixes $\gamma$ pointwise is
$\HH  = (U(n-1) \times \zun_{2})/  \zun_{n+1}$. We have an isomorphism $\HH \cong U(n-1)$.

After doing the above identifications 
$\HH \cong U(n-1)$, $\KK_1 \cong U(n-1)\times U(1)$ and 
$\KK_2 \cong U(n)$ the inclusion maps of $\HH$ in $\KK_1$ 
respectively $\KK_2$ are given by  
\[
A \mapsto (A,1) , \quad 
A \mapsto 
\begin{pmatrix}
A & 0\\
0 & 1\\
\end{pmatrix}
. 
\]

A simple closed geodesics $\gamma$ can be written as
\[ \gamma (t) = \rho ( \cos (\pi t) x + \sin (\pi t) v) , \quad 0 \leq t \leq1 \]
where $x\in S^{2n+1}$ and $v\in (\CC x)^\perp$  as described in \cite{GHL} example 2.110.
Furthermore, by this reference and lemma \ref{lemma:geodesics.symmetric}, the only points 
conjugate to $\gamma (0)$ is $\gamma (1/2)$ and $\gamma (1)$. We choose our coordinate
system such that $\gamma (0) = [e_{n+1}]$ and $\gamma (1/2) = [e_n]$.

For the geodesics $\gamma_m$, defined as the $m$-fold iterates of $\gamma$, 
the conjugate points of $\gamma_m (0)$ are $\gamma_m (i/2m)$, 
where $1\leq i\leq 2m$. 
But $\gamma_m (k/2m)=p$ if $k$ is even, and 
$\gamma_m (k/2m)=q$ if $k$ is odd.    
So the groups corresponding to $\gamma_m$ are given as follows:
\[
\KK (\gamma_m)_i=
\begin{cases}
\KK_1 &\text{if $i$ is odd}\\
\KK_2 &\text{if $i$ is even}\\
\end{cases}
\quad , \quad 1\leq i \leq 2m-1.
\]   
The associated $K$-cycle is 
\[
\kgcycle {\gamma_m} \GG = 
\EE (\GG , \KK_1 , \KK_2 , \KK_1 ,\dots , \KK_2 , \KK_1 ;\HH ).
\]
This manifold contains a copy of $\GG /\HH$ by the section which  
maps $g\HH$ to $[g,e, \dots ,e]$. This section is the inclusion 
of a submanifold of $\kgcycle {\gamma_m} \GG$ with normal bundle 
$\nu_m$. There is an associated Pontryagin-Thom collapse map
$c: {\kgcycle {\gamma_m} \GG }_+ \to \Th(\nu_m)$. 

We now verify the conditions of theorem \ref{thm:splitting} part 2. with $k=1$ suspension. 
The group $\KK_1$ acts transitively on the unit circle $S^1\subseteq \CC$ by $((A,z),u)\mapsto zu$.
The isotropy group at $1\in S^1$ is $U(n-1) \cong \HH$. So there
is an isomorphism $\KK_1/ \HH \cong S^1$. The left $\HH$-action corresponds to the
trivial $\HH$-action on $S^1$ under this identification. By Example \ref{ex.repr}, the fixed point
$[e]\in (\KK_1 /\HH )_+$ is splitting up to 1-fold suspension.

The group $\KK_2$ acts transitively on the unit sphere $S^{2n-1}$ in $\CC^n$ by the standard 
action $(A,x)\mapsto Ax$. The isotropy group at $e_n$ is $U(n-1) \cong \HH$. So there
is an isomorphism $\KK_2/ \HH \cong S^{2n-1}$. Under this identification, the left $\HH$-action on 
$S^{2n-1}$ is via the standard inclusion of $U(n-1)$ in $U(n)$. The fixed point $[e]\in (\KK_2 /\HH )_+$,
which corresponds to $e_n\in S^{2n-1}$, is splitting up to 1-fold suspension by Example \ref{ex.repr}.
 
Notice that $W_1=T_{[e]} (\KK_1 /\HH)$ is a trivial 1-dimensional
$\HH$-representation. The tangent space $T_{e_n}(S^{2n-1})$ decomposes as a direct sum
of the subspace $(\CC e_n)^\perp$ with complex structure and basis $e_1, \dots , e_{n-1}$ and a one dimensional
real subspace with basis $ie_n$. The $\HH$-action can be found via the curves 
$\gamma (t) = e_n \cos t + v \sin t$, with $v\in T_{e_n}(S^{2n-1})$. 
One sees that the $\HH$-representation
$W_2=T_{[e]} (\KK_2 /\HH )$ is the standard $2n-2$-dimensional $U(n-1)$-representation 
with a trivial 1-dimensional representation added.

The $\HH$-representations $W_1$ and $W_2\oplus \epsilon^1$ both extend
to representations of $\KK_2 \cong U(n)$. Namely, the trivial $1$-dimensional representation 
and the standard $2n$-dimensional representation. So, in order to check the remaining conditions 
of theorem \ref{thm:splitting} part 2. it is enough to see that the $\HH$-representations
\[
(W_1\oplus W_2)^{\oplus j} \oplus W_1 , \quad
W_2\oplus (W_1\oplus W_2)^{\oplus j} \oplus W_1, \quad j\geq 0
\]
extends to $\KK_2$ representations. But $W_1\oplus W_2$ extends to 
a representation of $\KK_2 \cong U(n)$, namely the standard representation so this is OK. 
By theorem \ref{thm:splitting} we now get the following result:

\begin{theorem}
\label{th.complex.splitting}
Let $p: S(\tau )\to \CP^n$ be the unit sphere bundle of the
tangent bundle $\tau$ on $\CP^n$. Let $\xi_m$ be the vector bundle 
$(p^*(\tau ))^{\oplus (m-1)} \oplus \epsilon$
over $S(\tau )$, where $\epsilon$ is a one dimensional trivial real vector bundle.
Then, there is a homotopy equivalence of spaces
\[
\Sigma (\Lambda \CP^n)_+ \simeq 
\Sigma \CP^n_+ \vee \bigvee_{m=1}^\infty \Sigma \Th (\xi_m). 
\]
\end{theorem}
\begin{proof}
We have to check that all non-trivial critical submanifolds have positive 
dimensional negative bundles. We know all critical points for the energy 
functional (i.e. the closed geodesics) by lemma \ref{lemma:geodesics.symmetric}. The index of one 
of these is given by equation (\ref{eq:dimension}). By the calculations above,
$\dim(\KK_1)-\dim(\HH)=1$, so it follows that
the index of a critical point which is not a
constant path is greater or equal to 1.

It only remains to identify the negative bundles $\norbd_\nu^-$
with the bundles $\xi_m$.
We have identified the negative bundle corresponding
to $\gamma_m$ with the bundle over $\GG /\HH$ induced by 
$V=(W_1\oplus W_2)^{\oplus (m-1)}\oplus W_1$.
So it suffices to see that the bundle $\eta $ over $\GG /\HH$ induced
by $W_1\oplus W_2$ agrees with $p^*(\tau )$.

By equation (\ref{eqn:tangent_space_id}) the maps $d\rho _x$ gives a diffeomorphism 
$V_2(\CC^{n+1} )/\zun_2 \cong S(\tau )$ where $V_2(\CC^{n+1} )$ is the Stiefel manifold
of orthonormal $2$-frames in $\CC^{n+1}$. The $\GG$-action $[A]\cdot [u,v]=[Au,Av]$ 
is transitive with isotropy group at $[e_{n+1}, e_n]$ equal to $\HH$. So there is an
isomorphism $\GG /\HH \cong S(\tau )$ such that the base spaces of the two bundles agree. 

The $\GG$-action on $\CP^n$ given by $[A]\cdot [u]=[Au]$ is transitive 
with isotropy group at $[e_{n+1}]$ equal to $\KK_2$. So there is also an isomorphism 
$\GG /\KK_2 \cong \CP^n$. 
Under these identifications, the projection $p$ corresponds to the standard identification map
$\GG /\HH \to \GG /\KK_2$, or equivalently to the map 
\[ 
U(n+1)/ (U(n-1)\times \zun (2)) \to U(n+1)/(U(n)\times U(1)).
\]
The bundle $\eta$ is isomorphic to the pullback under this map of the bundle $\overline \eta$
induced by the standard representation of $U(n)$. However, there is an isomorphism 
$\overline \eta \cong \tau$. The fibers of the tangent bundle are given by 
$T_{\rho(x)} (\CP^n )\cong (\CC x )^\perp$ 
and the isomorphism is given by the bundle map
\[
U(n+1) \times_{(U(n) \times U(1))} \CC^n \to \tau ; \quad
[A,v] \mapsto ([A e_{n+1}], A \begin{pmatrix}
v \\
0 
\end{pmatrix} ).
\]
This finishes the proof of the theorem.
\end{proof}

\subsection{The quaternionic projective spaces}

The standard metric on $\HP^n$ can be described in a similar way as the Fubini-Study metric
on $\CP^n$. For the quaternions $\mathbb{H}$ the conjugate is defined as
\[ z=x_0+x_1i+x_2j+x_3k \mapsto z^* = x_0-x_1i-x_2j-x_3k. \] 
It satisfies the rule $(ab)^* = b^* a^*$. The associated norm $|z|=\sqrt{zz^*}$ agrees with the norm of 
$z$ viewed as a 4-dimensional vector. 

Equip $\mathbb{H}^{n+1}$ with standard inner product 
$\langle v, w \rangle = \sum v_\ell w_\ell^*$. The group of matrices over $\mathbb{H}$ which 
preserves this structure is denoted $Sp(n+1)$. Let $S^{4n+3}$ be the
unit sphere in $\mathbb{H}^{n+1}$ and let
\[ \rho : S^{4n+3} \to S^{4n+3}/Sp(1) = \HP^n \]
be the canonical projection. The standard metric is determined by the property that
\[ d\rho_x : (\mathbb{H} x)^\perp \hookrightarrow T_x(S^{4n+3}) \xrightarrow{D_x (\rho )} T_{\rho (x)} (\HP^n ) \]
is an $\RR$-linear isometry for all $x\in S^{4n+3}$. 

The group $Sp(n+1)$ acts by isometries on $\HP^n$ and every isometry is induced by a matrix from this group.
The kernel of the action is however not the group $\{ zI| z\in Sp(1) \}$ but smaller. As 
$\langle zv , zw \rangle = z \langle v, w \rangle z^*$ one gets the further condition that $z$ must lie in the center of
$\mathbb{H}$ such that the kernel equals $\{ \pm I \}$. 
Thus the isometry group can be identified as $\GG = Sp(n+1)/\{\pm I\}$. 

The space $\HP^n$ is isotropic and as in the
complex case we get that the isotropy group of the
geodesic $\gamma$ starting at $[e_{n+1}]$ in the direction $e_n$ is
\[
\HH = Sp(n-1) \times Sp(1)/\{\pm I\}.
\]
The isotropy group of the point $[e_{n+1}]$ is
\[
\KK_2=Sp(n) \times Sp(1)/\{\pm I\},
\]
and the isotropy group of the pair of conjugated points
$[e_n],[e_{n+1}]$ is
\[
\KK_1=Sp(n-1)\times Sp(1) \times Sp(1)/\{\pm I\}.
\]
The inclusion maps $\HH\subset \KK_1$ and $\HH\subset \KK_2$ 
are given as follows:
\[
(A_1,A_2) \mapsto (A_1,A_2,A_2) , \quad 
(A_1,A_2)\mapsto (
\begin{pmatrix}
A_1 & 0\\
0 & I \\
\end{pmatrix}
, A_2) 
\]
Any closed geodesic will be obtained from a simple closed 
geodesic by running through it $m$ times. 
As in the complex case, the isotropy groups associated to 
the conjugated points will alternate between
$\KK_1$ and $\KK_2$. 

The quotient spaces 
\begin{align*}
& \KK_1 /\HH \cong 
Sp(n-1)\times Sp(1)\times Sp(1)/ (Sp(n-1)\times Sp(1)) \cong Sp(1) ,\\
& \KK_2 /\HH \cong 
Sp(n)\times Sp(1) /(Sp(n-1)\times Sp(1))\cong Sp(n)/Sp(n-1)
\cong S^{4n-1}
\end{align*}
are equivariantly unit spheres in $\HH$-representations. So by 
Example \ref{ex.repr} the fixed points $[e]\in (\KK_1 /\HH)_+$ and
$[e]\in (\KK_2 /\HH)_+$ are splitting up to $1$-fold suspension. We look closer at the 
representations.

The $\HH$-representation 
$W_1=T_{[e]}(\KK_1 /\HH )\cong \mathfrak{sp}(1)$
is induced from the adjoint representation of $Sp(1)$ under
the projection map 
$$\HH =Sp(n-1)\times Sp(1)/\{\pm I\} \xrightarrow{p_2} Sp(1)/\{\pm I\}.$$
The adjont representation of $Sp(1)$ can be described as the $3$-dimensional vector space of purely
imaginary quaternions $\tilde \HHQ$ on which $Sp(1)$ acts by conjugation.

The $\HH$ representation $W_2=T_{[e]}(\KK_2 /\HH )
\cong \mathfrak{sp}(n) / \mathfrak{sp}(n-1)$ is equivalent to
the representation of $\HH \cong Sp(n-1)\times Sp(1)/\{\pm I\}$ on
$\HHQ^{n-1} \oplus \tilde \HHQ$, given by
\[
(A, B)(v, w)=(AvB^{-1}, BwB^{-1}).
\]
This representation splits as a sum of two representation, one of 
them being given by the adjoint representation of the $Sp(1)$ 
factor in $\HH$.

The representation $W_1$ extends to a 
$\KK_1\cong Sp(n-1)\times Sp(1)\times Sp(1)/\{\pm I\}$ and to a
$\KK_2\cong Sp(n)\times Sp(1)/\{\pm I\}$ representation, by projecting the groups onto
the last factor $Sp(1)$, and composing with the adjoint
representation.

The representation $W_2$ does not itself extend, but if you add
a trivial 1-dimensional representation to it, it does.
Put 
\[
\HHQ^n =\HHQ^{n-1} \oplus \HHQ =
\HHQ^{n-1} \oplus \tilde \HHQ \oplus \epsilon \cong 
W_2\oplus \epsilon.
\] 

\begin{theorem}
\label{th.quaternionic.splitting}
Let $p: S(\tau )\to \HP^n$ be the unit sphere bundle of the
tangent bundle $\tau$ on $\HP^n$. Let $\eta$ be the vector
bundle over $\HP^n \cong S^{4n+3} /Sp(1)$ induced by the conjugation action of $Sp(1)$
on the $3$-dimensional real vector space of purely imaginary quaternions $\tilde \HHQ$.
Let $\xi_m$ be the vector bundle 
$p^*(\tau^{\oplus (m-1)} \oplus \eta^{\oplus m})$ over $S(\tau )$. 
Then there is a homotopy equivalence of spectra
\[
\Sigma^\infty (\Lambda \HP^n)_+ \cong \Sigma^\infty \HP^n_+
\vee
\bigvee_{m\geq 1} \Sigma^{\infty-m+1} \Th(\xi_m ).
\]
\end{theorem}

\begin{proof}
We easily check the conditions of
theorem \ref{thm:splitting} (compare to theorem
\ref{th.complex.splitting})
so we only have to see that
there is a bundle equivalence between $\xi_m$
and the pullback of the bundle induced by 
$(W_2\oplus \epsilon)^{\oplus (m-1)} \oplus W_1^{\oplus m}$.

Both these bundles are induced from bundles on $\HP^n$, so it
is sufficient to check the isomorphism on the appropriate bundles there.

But the bundle $W_1$ is by definition $\eta$, so we only need to show
that the bundle induced by $W_2\oplus\epsilon$ on $\HP^n$ is
the tangent bundle. We already identified $W_2\oplus \epsilon$ with
the standard representation $\HHQ^n$. And now we can give an
isomorphism as a bundle map
\[
Sp(n+1) \times_{Sp(n)\times Sp(1)} \HHQ^n \to \tau \quad ; \quad
[A,v] \mapsto ([A e_{n+1}], A \begin{pmatrix}
v \\
0 
\end{pmatrix} ).
\]
\end{proof}

\subsection{The Cayley projective plane}
\label{sec:Cayley}

Now we consider the Cayley projective plane $\CaP$. Recall from
\ref{Section:Cayley} that $\CaP$ is a symmetric space of rank 1.
The isometry group $F_4$ acts isotropically. The isotropy
group of a point $x$ is $\KK =Spin(9)$, the isotropy group
of a geodesic $\gamma$ is denoted by $\HH$. 
We saw in theorem \ref{th:Cayley.geodesics} that the isotropy groups
of pairs of conjugate points $(\gamma (0),\gamma (t))$ on $\gamma$ 
alternate between
$\KK_1 =Spin(8)$ and $\KK_2 \cong \KK=Spin(9)$.
Let $W_1$ and $W_2$ denote the $\HH$-representations 
$T_{[e]}(\KK_1 /\HH)$ and $T_{[e]}(\KK_2 /\HH)$ respectively. 
Note that $[e]\in \KK_i /\HH$ is not a fixed point under the
$\KK_i$-action, so a priory we don't even know that $W_i$
extends to a $\KK_i$-representation.

\begin{lemma}
\label{lemma:Cayley.representation.theory}
The action of $\HH$ on $W_2\oplus \epsilon$ can be extended to an 
action of $\KK_2$, and thus to an action of the subgroup $\KK_1$.
The class $[W_1]\in RO(\HH)$ is in the image of the restriction map
$RO(\KK)\to RO(\HH)$, and thus in the image of $RO(\KK_1)\to RO(\HH)$.
\end{lemma}

\begin{proof}
Even if the action of $\HH$ on  $W_2$ does not itself extend to an
action of $\KK_2$, it does in the stable sense. To see this, note that
$\KK_2$ acts transitively on $S(T_x(\CaP ))$ with isotropy group $\HH$ at the 
unit vector corresponding to $\gamma$. 
Thus, $\KK_2 / \HH \cong S(T_x(\CaP ))$. This means that  
the sum $W_2\oplus \epsilon$ is the tangent space $\HH$-representation $T_x(\CaP )$, which by
definition is the restriction of a representation of $\KK_2$.

The chain of inclusions of Lie groups $\HH \subset \KK_1 \subset \KK_2$
defines a short exact sequence of $\HH$-representations.
\[
0\to T_{[e]}(\KK_1 /\HH ) \to T_{[e]}(\KK_2 /\HH ) \to 
T_{[e]}(\KK_2 /\KK_1 ) \to 0
\]
It follows that $W_2\cong W_1 \oplus T_{[e]}(\KK_2 /\KK_1 )$ as 
$\HH$-representations. So we have to prove that 
\[
[T_{[e]}(\KK_2 /\KK_1 )]\in \mathrm{Im} (RO(\KK_2 )\to RO(\HH )).
\]
By definition, $[T_{[e]}(\KK_2 /\KK_1 )]$ is the restriction 
of a $\KK_1$-representation, so we have to show that
\[
[T_{[e]}(\KK_2 /\KK_1 )]\in \mathrm{Im} (RO(\KK_2 )\to RO(\KK_1 )).
\]
However, according to \ref{th:Cayley.geodesics}, we can identify 
$\KK_1$ with $Spin(8)$ and $\KK$ with $Spin(9)$ in such a way that 
the inclusion $\KK_1 \subset \KK$ corresponds to the standard
inclusion $Spin(8)\subset Spin(9)$. Under this isomorphism, the 
representation $T_{[e]}(\KK_2 /\KK_1 )$ corresponds to 
$T_{[e]}(Spin(9)/Spin(8))$ which is just the standard 8-dimensional 
representation $\rho_8$ of $Spin(8)$. The lemma follows from
the fact that the restriction of the standard $Spin(9)$-representation 
$\rho_9$ to $Spin(8)$ is $\rho_8\oplus \epsilon$.
\end{proof}
 
\begin{remark}
  In the proof of lemma
  \ref{lemma:Cayley.representation.theory} we intentionally 
avoid the discussion of precisely which $\HH=Spin(7)$
representations we are dealing with. This would involve
a treatment of triality, which we do not want to
bring up here. 
\end{remark}

\begin{theorem}
\label{th:cayley.splitting}
Let $p:S(\tau )\to \CaP$ be the unit sphere bundle of the tangent 
bundle $\tau$ on $\CaP$. Let $\eta$ be the 7-dimensional vector bundle 
induced from $W_1$ over $\CaP \cong F_4/Spin(9)$. Let $\norbd_m$ be the 
vector bundle $p^*(\tau^{\oplus (m-1)}\oplus \eta^{\oplus m})$ over 
$S(\tau )$. Then there is a homotopy equivalence of spectra
\[
\Sigma^\infty (\Lambda \CaP_+)\simeq \Sigma^\infty \CaP_+ 
\vee
\bigvee_{m=1}^\infty \Sigma^{\infty -m+1}\Th (\norbd_m ).
\]  
\end{theorem}

\begin{proof}
This is analogous the proof of theorem
\ref{th.quaternionic.splitting}. We have to check that
the bundle induced by $W_2\oplus \epsilon$ agrees with the 
tangent bundle. It is sufficient to check that the
$\KK_2$ representations on the tangent space $T_{\gamma(0)}(\CaP )$
agrees with the $\KK_2$-representation $W_2\oplus \epsilon$.
But actually, it is well known that the tangent representation of 
$Spin(9)$ acting as the isometry group fixing a point on $\CaP$ 
is exactly the sum of a trivial representation and the spinor 
representation of $Spin(9)$.
\end{proof}

\section{Comparison with earlier results}
\label{sec:Thom}

We want to compare the results of this paper with the results of
\cite{BO}, which motivated it. 
We first recall some notation used in \cite{BO}.

Let $X$ be either $\CP^n$, $\HP^n$ or $\CaP$. Let $\tau$ be 
the tangent bundle over $X$ and let 
$q\tau= \tau \oplus \dots \oplus \tau$ denote the $q$-fold 
Whitney sum. We define $C_q(X)$ to be the cofiber of the cofibration
\[
\Th (q\tau ) \to \Th ((q+1)\tau ).
 \]
The cohomology of $X$ is a truncated polynomial ring. 
The generator has degree $r(X)$, where $r(X)=2,4,8$ respectively.
In \cite{BO}, we proved the following theorem:
\begin{theorem}
\label{th.BO}
There is an isomorphism of modules over the mod two Steenrod algebra
\[
H^*(\Lambda X; \ZZ /2) \cong 
H^*(X_+\vee \bigvee_{q\geq 0} \Sigma^{(r(X)-2)(q+1)} C_q(X);\ZZ /2).
\]
\end{theorem}

To reinterpret this in terms of the bundles considered in this paper,
we need the following rewriting:
Let $X$ be a space and $\xi_1 , \xi_2$ vector bundles over $X$. 
Let $D(\xi_i )\to X$, $S(\xi_i )\to X$, $\Th (\xi_i )$ be the 
corresponding disk bundle, sphere bundle, and Thom space.
Let $\epsilon$ denote a trivial line bundle.

\begin{theorem}
\label{th.rewriting}
Let $C(\xi_1 ,\xi_2 )$ be the cofiber of the map 
$\Th (\xi_1 ) \to \Th (\xi_1 \oplus \xi_2 )$
given by inclusion of $\xi_1$ and the zero section of $\xi_2$.
Let $p:S(\xi_2 )\to X$ be the projection map for the sphere bundle 
of $\xi_2$. Then there is a homotopy equivalence
$$C(\xi_1 , \xi_2 ) \simeq \Th (p^*(\xi_1 \oplus \epsilon )).$$
\end{theorem}

\begin{proof}
There is a diagram of cofibrations
\[
\begin{CD} 
S(\xi_1 )_+ @>>> D(\xi_1 )_+ @>>> \Th (\xi_1 ) \\
@VVV @VVV @VVV \\
S(\xi_1 \oplus \xi_2 )_+  @>>> D(\xi_1 \oplus \xi_2 )_+
@>>> \Th (\xi_1 \oplus \xi_2 ) \\
@VVV @VVV @VVV \\
S(\xi_1 \oplus \xi_2 )/S(\xi_1 ) @>>> D(\xi_1 \oplus \xi_2 )/D(\xi_1 )
@>>> C(\xi_1 , \xi_2 )
\end{CD}
\]

The inclusion $D(\xi_1) \to D(\xi_1 \oplus \xi_2)$
is a homotopy equivalence, so its cofiber is contractible. 
Thus by extending the above diagram to the right, we get a 
homotopy equivalence 
$C(\xi_1 , \xi_2 ) \simeq \Sigma S(\xi_1 \oplus \xi_2 )/S(\xi_1 )$.
Next, consider the inclusion 
$$S(\xi_1 ) \hookrightarrow S(\xi_1 \oplus \xi_2 ) \simeq 
S(\xi_1)\times_X D(\xi_2 ) \cup D(\xi_1 ) \times_X S(\xi_2 ).$$ 
Since $S(\xi_1 ) \hookrightarrow S(\xi_1 ) \times_X D(\xi_2 )$ 
is a homotopy equivalence we have that the cofiber of the map
$f:S(\xi_1 ) \times_X S(\xi_2) \to D(\xi_1 )\times_X S(\xi_2 )$
is homotopy equivalent to $S(\xi_1 \oplus \xi_2 )/S(\xi_1 )$.

But $S(p^* (\xi_1 )) \simeq S(\xi_1 )\times_X S(\xi_2 )$
and $D(p^* (\xi_2 )) \simeq D(\xi_1 )\times_X S(\xi_2 )$
such that $\Th (p^* (\xi_1 ))$ is the cofiber of $f$.
The result follows.
\end{proof}

We use this result to rewrite theorem \ref{th.BO}.
\begin{theorem}
\label{th.BO1}
The cohomology $H^*(\Lambda X ; \ZZ/2)$ is as module over the 
mod two Steenrod algebra isomorphic to the spectrum cohomology 
\[
H^*(X_+\vee \bigvee_{m\geq 1} 
\Sigma^{-m+1} \Th (p^*((m-1)\tau \oplus m\epsilon^{r(X)-1} );\ZZ /2).
\]
\end{theorem}

\begin{proof}
This follows from applying theorem \ref{th.rewriting} to
$\xi_1 =q\tau $ and $\xi_2 =\tau$, rewriting the Thom space a little,
and finally substituting $m=q+1$.
\end{proof}

We now check that this agrees with the results obtained in 
section \ref{sec:Rank1}.
\begin{remark}
The case $X=\CP^n$, $r(X)=2$, follows directly from theorem
\ref{th.complex.splitting} after identifying the bundle $\xi_m$.   
\end{remark}

\begin{remark}
\label{rem.quaternionic}
In the case $X=\HP^n$, $r(X)=4$ theorem
\ref{th.quaternionic.splitting} says that if $\eta$ is the 
bundle induced from the adjoint representation of $S^3$ on
$\HP^n$, then $H^*(\Lambda \HP^n ;\ZZ /2)$ is isomorphic as a
module over the Steenrod algebra to
\[
H^*(X;\ZZ /2)\oplus \bigoplus_{m\geq 1} 
\tilde H^*(\Sigma^{-m+1} \Th (p^*((m-1)\tau \oplus m\eta );\ZZ /2).
\] 
So it is sufficient to see that there is an isomorphism of
modules over the Steenrod algebra
\[
H^*(\Th (p^*((m-1)\tau ) \oplus mp^*(\eta )) ;\ZZ /2) \cong
H^*(\Th (p^*((m-1)\tau ) \oplus m\epsilon^3 );\ZZ /2).
\]
According to \cite{Characteristic}, is is sufficient to see 
that the Stiefel Whitney classes of $p^*(\eta )$ vanish. By
naturality, it is sufficient to see that the
Stiefel Whitney classes 
\[
w_j(\eta )\in H^*(\HP^n ;\ZZ /2)
\]
vanish. But if $j\leq 3$, the group
$H^j(\HP^n )$ vanishes, so $w_j (\eta )=0$. If 
$j>3$, the class $w_j(\eta)$ vanish, since the bundle 
$\eta$ is 3 dimensional.  
\end{remark}

\begin{remark}
As in remark \ref{rem.quaternionic} we see from 
theorem \ref{th:cayley.splitting} that in order
to check the part of theorem \ref{th.BO1} involving
the Cayley projective plane, it is sufficient
to see, that the Stiefel Whitney classes of the bundle
$\eta$ over $S(\tau )$ are trivial. The bundle is 7-dimensional 
so $w_j(\eta )=0$ for $j\geq 8$.
On the other hand, $w_j\in H^j( S(\tau );\ZZ /2)=0$
for $1\leq j \leq 7$, so all Stiefel Whitney classes are
indeed trivial.    
\end{remark}

\end{document}